\documentclass[11pt]{article}
\usepackage{amsmath}
\usepackage{amssymb,amsmath}
\usepackage{mathrsfs}
\usepackage{amsmath,amsthm,amsfonts,amssymb,bm,euscript}
\usepackage{fancyhdr,amscd}
\usepackage{anysize}
\usepackage{graphicx,epsfig}
\usepackage{amsthm,latexsym,amsmath,amssymb}
\usepackage{amssymb,amsmath}
\usepackage{mathrsfs}
\usepackage{graphicx}
\usepackage{color}
\usepackage{mathrsfs}
\usepackage{cite}
\usepackage{indentfirst}
 \usepackage[titletoc]{appendix}
 \usepackage [latin1]{inputenc}
\usepackage{amsthm, amsfonts, mathrsfs}
\usepackage{mathptmx}
\usepackage{fullpage}
\usepackage{amsfonts,graphicx}
\usepackage[colorlinks=true,  linkcolor=blue, citecolor=blue, urlcolor=blue]{hyperref}
\let\al=\alpha

\let\p=\psi

\def\na{\nabla}

\def\p{\partial}

\def\eqdefa{\buildrel\hbox{\footnotesize def}\over =}

\def\C{\mathop{\bf C\kern 0pt}\nolimits}
\def\DD{\mathop{\bf D\kern 0pt}\nolimits}
\def\K{\mathop{\bf K\kern 0pt}\nolimits}
\def\N{\mathop{\bf N\kern 0pt}\nolimits}
\def\Q{\mathop{\bf Q\kern 0pt}\nolimits}
\def\R{\mathop{\bf R\kern 0pt}\nolimits}
\def\endproof{\hphantom{MM}\hfill\llap{$\square$}\goodbreak}

\newcommand{\Dv}{{\rm div}}

\newcommand{\beq}{\begin{equation}}
\newcommand{\eeq}{\end{equation}}
\newcommand{\ben}{\begin{eqnarray}}
\newcommand{\een}{\end{eqnarray}}
\newcommand{\beno}{\begin{eqnarray*}}
\newcommand{\eeno}{\end{eqnarray*}}

\newtheorem{Theorem}{Theorem}[section]
\newtheorem{Definition}[Theorem]{Definition}
\newtheorem{Proposition}[Theorem]{Proposition}
\newtheorem{Lemma}[Theorem]{Lemma}

\newtheorem{Remark}[Theorem]{Remark}

\numberwithin{equation}{section}

\allowdisplaybreaks \numberwithin{equation} {section}

\begin{document}
\title{Fujita-Kato solution for  the 3D  compressible pressureless Navier-Stokes  equations with  discontinuous and large-variation density
}
\author{ Xiaojie  Wang$^{1}$ \ \  Jiahong Wu$^{2}$ \ \   Fuyi  Xu$^{1\dag}$\\[2mm]
 { \small   $^1$School of Mathematics and  Statistics, Shandong University of Technology,}\\
  { \small Zibo,    255049,  Shandong,    China}\\
  { \small  $^2$Department of Mathematics, University of Notre Dame,}\\
   { \small Notre Dame, IN 46556, USA}
   }
         \date{}
         \maketitle
\noindent{\bf Abstract}\ \ \ This paper mainly  focuses on the Cauchy problem to the 3D  compressible pressureless Navier-Stokes  equations arising from
models of collective behavior, which can be derived by taking the high Mach number limit of the classical  compressible Navier-Stokes system. We construct  the  global-in-time existence and uniqueness of the so-called Fujita-Kato solution to the system,  provided that the initial density $\rho_0$ is discontinuous, large-variation and the initial velocity $u_0$ is in a critical functional framework. Our method relies on some time weighted estimates  and the Lagrangian approach. 

\vskip   0.2cm \footnotetext[1]{$^\dag$Corresponding author.}
\vskip   0.2cm \footnotetext[2]{E-mail addresses: wxj15615637750@163.com(X. Wang), \ \  jwu29@nd.edu(J. Wu), \ \ zbxufuyi@163.com(F. Xu).} \setlength{\baselineskip}{20pt}

\section{Introduction and Main Results}
\setcounter{section}{1}\setcounter{equation}{0}  \subsection{Model and synopsis of related studies} \ \
\ \  In recent years, due to the physical importance and the very broad applications of
hydrodynamics in nature and industry,   the following general form of hydrodynamical
system has received a extensive attention
in the mathematical community
\begin{align}\label{1.1}
\left\{
\begin{aligned}
&\partial_t \rho+\nabla\cdot(\rho u)=0,\\
&\partial_t (\rho u)+\nabla\cdot(\rho u\otimes u)+\nabla P(\rho)=\mathcal{D}(u,\rho).
\end{aligned}
\right.
\end{align}
Here $\rho=\rho(t,x)$ and $u=u(t,x)$ are the unknown functions representing the density and velocity, respectively,  $P(\rho)$ stands for the pressure.
Hydrodynamic limiting system \eqref{1.1} with different type of pressures including the pressureless dynamics with $P\equiv0$ \cite {FK} and isothermal pressure with  $P(\rho)=\kappa \rho^{\gamma} (\kappa>0, \gamma>1)$  \cite{KMT},  can be formally
derived from  the following system
\begin{equation}\label{micro1}
  \partial_t f+v\cdot \nabla_x f+\nabla_v\cdot K(f)f=0,
\end{equation}
where $f = f (t, x, v)$ is a distribution function of a gas in the phase-space. Classically,
if the operator $K$ comes from the Poisson potential, then we find the Vlasov system.
If taking a less singular operator, then one may obtain for example the Cucker-Smale
system that models collective behavior like flocking of birds.

When the term $\mathcal{D}(u,\rho)$ represents the nonlocal velocity alignment which is given as follows
\begin{equation}\label{pressure1}\mathcal{D}(u,\rho)=-\int_{\mathbb{R}^{N}}\phi(x-y)(u(x)-u(y))\rho(y,t)\rho(x,t)\mathrm{d}y,
\end{equation}
where $\phi$ is called the communication weight, measuring the
strength of the alignment interactions. Then the system \eqref{1.1} is called the  the Euler-alignment equations. An interesting type of communication weights are strongly singular near the origin,
with a prototype taking the following form
\begin{equation}\label{phi}
  \phi(x)=\frac{c_{\alpha,N}}{|x|^{N+\alpha}},\quad
   c_\alpha=\frac{2^{\alpha}\Gamma(\frac{\alpha+N}{2})}{\pi^{\frac{N}{2}}|\Gamma(-\frac{\alpha}{2})|}.
\end{equation}
where the constant $c_\alpha$ is related to the fractional Laplacian
operator, which in $\mathbb{R}^N$ reads
\begin{equation}\label{LD}
  \Lambda^{\alpha}=(-\Delta)^{\frac{\alpha}{2}},\quad
  \Lambda^{\alpha}f=c_{\alpha,N}P.V.
  \int_{{\mathbb R}^N}\frac{f(x)-f(y)}{|x-y|^{N+\alpha}}\mathrm{d}y,
  \quad 0<\alpha<2.
\end{equation}
It is evident that  the singular alignment $\mathcal{D}(u,\rho)$ with such a weight $\phi$ can be
expressed as a commutator form related to the fractional Laplace operator $\Lambda^{\alpha}$
\begin{equation}\label{pressure2}\mathcal{D}(u,\rho)=-\rho(\Lambda^{\alpha}(\rho u)-u\Lambda^{\alpha}\rho)=-\rho([\Lambda^{\alpha},u]\rho).
\end{equation}
The  pressureless Euler-alignment equations (i.e., $P\equiv0$)  has an intriguing regularization effect to the solutions. Global regularity
is obtained for any non-viscous smooth initial data for the system in the one-dimensional
torus by Shvydkoy and Tadmor \cite{ST} for $1\leq \alpha< 2$, and by Do \emph{et al} \cite{DKRT} for $0 < \alpha < 1$ (see also
[40]). For the multi-dimensional case, global well-posedness
are only known for small initial data around an equilibrium state,  and  we refer readers to Refs \cite{DMPW,Sh}  and references therein.  Global
regularity for general large initial data remains a challenging open problem.

 When $\mathcal{D}(u,\rho)=0$,  the system \eqref{1.1} reduces to  the compressible Euler system, only local
existence of smooth solutions could be expected, and shock waves will form in finite time. If
$P\equiv0$, then one just recovers the pressureless compressible Euler system, and there is
no interaction whatsoever between the individuals, which describe  the formation of large-scale structures in astrophysics and the aggregation
of sticky particles  (see, e.g., \cite{SSZ, Z, YZ}) and the references therein. The majority of results in the literature (e.g., \cite{Bou1,Bou2,Bre1,Bre2,HW,Hynd}) focus on the mathematical theory of the   pressureless dynamics.

When  we replace the dissipation term $\mathcal{D}(u,\rho)$ by the
Lam\'{e} operator $-\mu\Delta u-(\lambda+\mu)\nabla\Dv u $, where  the constant viscosity coefficients $\mu$ and $\lambda$ comply with the physical conditions $\mu>0$  and $2\mu+N\lambda\geq0$,  then the system \eqref{1.1} becomes the following  classical barotropic compressible Navier-Stokes equations
\begin{align}\label{1.11}
\left\{
\begin{aligned}
&\p_{t}\rho+\Dv (\rho u)=0,\\
&\p_{t}(\rho u)+{\Dv}(\rho u\otimes u)-\mu\Delta u-(\lambda+\mu)\nabla\Dv u+\nabla P(\rho)=0.
\end{aligned}
\right.
\end{align}
The model \eqref{1.11} is  a fundamental mathematical model used to describe the motion of a viscous flow. Due to its mathematical challenges and broad physical applications, there are extensive literature on the well-posedness of solutions  for this system.
 Here we briefly review some of the most
relevant papers about global well-posedness of strong solutions to the system when the initial data are close to a constant equilibrium state.
Matsumura and Nishida \cite{MN1,MN2}  studied the global existence of classical
solutions to the 3D compressible Navier-Stokes equations   for data $(\rho_{0}, u_{0})$ with high regularity order and close to a stable
equilibrium. Motivated by those works on the
incompressible Navier-Stokes system \cite{FK1},  Danchin \cite{Dan1}  proved   the existence  and uniqueness of the global strong solution for  initial data  close to a stable
equilibrium state  in $L^{2}$ critical Besov spaces. Subsequently, the result of \cite{Dan1} was independently extended
to more general $L^{p}$ critical Besov spaces by Charve and Danchin \cite{CD},   Chen et al. \cite{CMZ2} and Haspot \cite{HS2}. It should be noted that these theory requires the solution to have small oscillations from a uniform nonvacuum state, ensuring that the density is strictly away from the vacuum.

If we further consider  the high Mach number limit in  the system \eqref{1.11}, a pressureless viscous model then emerges.  More precisely, we first introduce the dimensionless Mach number defined as the ratio of the velocity $u$ to the reference sound speed. Then the rescaled triplet
\begin{equation*}
   \rho_{\varepsilon}(t,x)= \rho(\frac{t}{\varepsilon},x),\quad u_{\varepsilon}(t,x)=\frac{1}{\varepsilon}u(\frac{t}{\varepsilon},x),\quad \mu_{\varepsilon}= \frac{1}{\varepsilon}\mu,\quad \lambda_{\varepsilon}= \frac{1}{\varepsilon}\lambda
\end{equation*}
satisfies
\begin{align*}
\left\{
\begin{aligned}
&\p_{t}\rho_{\varepsilon}+\Dv (\rho_{\varepsilon} u_{\varepsilon})=0,\\
&\p_{t}(\rho_{\varepsilon} u_{\varepsilon})+{\Dv}(\rho_{\varepsilon} u_{\varepsilon}\otimes u_{\varepsilon})-\mu_{\varepsilon}\Delta u_{\varepsilon}-(\lambda_{\varepsilon}+\mu_{\varepsilon})\nabla\Dv u_{\varepsilon}+\nabla P(\rho_{\varepsilon})=0.
\end{aligned}
\right.
\end{align*}
Concerning that the Mach number $\varepsilon$ goes to infinity and still defining
\begin{equation*}
   \lim_{\varepsilon\to \infty}(\rho_{\varepsilon}, u_{\varepsilon},\mu_{\varepsilon},\lambda_{\varepsilon})=(\rho, u,\mu,\lambda),
\end{equation*}
we formally derive the following pressureless Navier-Stokes equations
\begin{align} \label{1.12}
\left\{
\begin{aligned}
&\rho_{t}+\Dv (\rho u)=0,\\
&\p_{t}(\rho u)+{\Dv}(\rho u\otimes u)-\mu\Delta u-(\lambda+\mu)\nabla\Dv u=0.
\end{aligned}
\right.
\end{align}
The system \eqref{1.12} does not
take into account the effects of internal pressure. A simple example is the motion of dust, or free particles, moving through space, as seen in astrophysics \cite{DD} or in multi-fluid systems \cite{Ber1}. Beyond inanimate matter, we can also consider models that
describe collective behavior, where particles--or rather agents--display some level of intelligence. In such cases, having a force like internal pressure is less natural, leading to the pressureless gas dynamics system. In present work, we mainly  focus on the system \eqref{1.12} in $\mathbb{R}^{+}\times \mathbb{R}^3$,  subject to the given initial data
\begin{equation}\label{1.13}
(\rho,u)|_{t=0}=(\rho_{0},u_{0}).
\end{equation}

In recent years, the study of  mathematical properties for the pressureless Navier-Stokes equations remains an important topic  because  the model has widespread applications in hydrodynamics, both in nature and industry \cite{Ber2, Bou, DD, ST, TT}.   
In 2021, Danchin et al. \cite{DMP}  first  proved global well-posedness of strong solution to the Cauchy problem \eqref{1.12}-\eqref{1.13} relying heavily on a set of maximal regularity estimates for the heat equation when the initial density $\rho_0$ is a small fluctuation of a positive constant in the $L^{\infty}(\mathbb{R}^{3})$-norm
and the initial velocity $u_0$  belongs to the  homogeneous Besov spaces (not critical regularity) satisfying a  smallness condition. This result underscores that the Cauchy problem \eqref{1.12}-\eqref{1.13} still has global existence and uniqueness, despite the inability to leverage the effective viscous flux. More recently, Guo  et al. \cite{GTZ} constructed the  global well-posedness and stability of classical solutions to the Cauchy problem \eqref{1.12}-\eqref{1.13} by virtue of the smallness of initial perturbation in Sobolev norm and the $L^{1}$ norm of initial velocity.  It is worthwhile to emphasize that,
 the assumption of small perturbations on the initial density plays a key role in their analysis, which excludes the case with large-variations along a hypersurface.
 Xu \cite{Xu} conducted a systematic study of the Gaussian bounds for the fundamental matrix of a generalized parabolic Lam\'{e} system with only bounded and measurable coefficients, subsequently establishing maximal $L^{1}$-regularity for the system.  As an application, he also proved global well-posedness for the Cauchy problem \eqref{1.12}-\eqref{1.13} when the  initial density  $\rho_0$  is large in  multiplier space $\mathcal{M}(\dot{B}_{p,1}^{\frac{3}{p}-1})$ and initial velocity $u_0$ is sufficiently small in critical Besov spaces $\dot{B}_{p,1}^{\frac{3}{p}-1}$.  Although Xu's approach accommodates initial density with large variations, it necessitates the use of multiplier space to require slightly more regularity on the initial density.

Based on the results mentioned above \cite{DMP,GTZ,Xu}, this raises the question of whether a small initial velocity in some critical space, combined with a large-variation, discontinuous bounded initial density, can yield a unique global-in-time solution to the 3D Cauchy problem \eqref{1.12}-\eqref{1.13}. The main contribution of the present paper is to provide a positive answer to this question.
\subsection{Main results} \ \
\ \  Our main result regarding global well-posedness is stated as follows.
\begin{Theorem}\label{1-1}
Let $(\rho_0,u_0)$ satisfy
\begin{equation}\label{1.2}
0<c_0\leq\rho_0(x)\leq C_0<+\infty\quad\text{and}\quad u_0\in \dot{B}_{2,1}^{\frac{1}{2}}.
\end{equation}
Then there
exists a small enough  constant $\varepsilon_0>0$ depending only on  $c_0, C_0$ such
that if
\begin{equation}\label{1.3}
\|u_0\|_{\dot{B}_{2,1}^{\frac{1}{2}}}\leq\varepsilon_0,
\end{equation}
the Cauchy problem \eqref{1.12}-\eqref{1.13} has a unique global  solution  $(\rho, u)$ with
$\rho\in L^\infty(\mathbb{R}^+\times\mathbb{R}^3)$ and $u\in C([0,\infty); \dot{B}_{2,1}^{\frac{1}{2}})\cap L^2(\mathbb{R}^+;\dot{B}_{2,1}^{\frac{3}{2}})$ which satisfies
\begin{equation}\label{1.4}
0<m\leq\rho(t,x)\leq M <\infty \quad \text{for}\quad(t,x)\in \mathbb{R}^+\times\mathbb{R}^3,
\end{equation}
and
\begin{equation}\label{1.5}
\begin{split}
&\|u\|_{\tilde{L}_{t}^{\infty}(\dot{B}_{2,1}^{\frac{1}{2}})}+\|u\|_{\tilde{L}_{t}^{2}(\dot{B}_{2,1}^{\frac{3}{2}})}+\|\sqrt{t}u\|_{\tilde{L}_{t}^{\infty}(\dot{B}_{2,1}^{\frac{3}{2}})}+\|\sqrt{t}\nabla u\|_{\tilde{L}_{t}^{2}(\dot{B}_{6,1}^{\frac{1}{2}})}+\|\sqrt{t}u_{t}\|_{\tilde{L}_{t}^{2}(\dot{B}_{2,1}^{\frac{1}{2}})}\\
&\quad+\|\sqrt{t}\nabla^{2}u\|_{L^{2}_{t}L^{3}}+\|tu_{t}\|_{\tilde{L}_{t}^{\infty}(\dot{B}_{2,1}^{\frac{1}{2}})}+\|tD_{t}u\|_{\tilde{L}_{t}^{2}(\dot{B}_{2,1}^{\frac{3}{2}})}\\
&\leq\|u_{0}\|_{\dot{B}_{2,1}^{\frac{1}{2}}}.
\end{split}
\end{equation}
In this context and throughout this paper, we denote  $D_t\eqdefa \p_t+u\cdot\na$ as the material derivative.
$m$ and $M $ represent two positive constants that depend on  $c_0$ and  $C_0$, respectively.
\end{Theorem}

 \begin{Remark}\label{1.1AA}
In contrast to the results in \cite{DMP,GTZ}, we eliminate the small perturbation assumption on the initial density, allowing for large variations in density. Meanwhile, the initial velocity $u_0$ has critical regularity.
\end{Remark}

 \begin{Remark}\label{1.1BB}
In this paper, we only require the density to be bounded, without imposing any regularity assumptions. As a result, our findings in Theorem \ref{1-1} expand upon the conclusions established in \cite{Xu}.
\end{Remark}

\begin{Remark}\label{1.1CC} Compared with  the existing results in \cite{QCZ,XQF,P}, we here deal with the compressible flow that the divergence is no longer 0,
 and need further exploit the  key  Lipschitz estimate of the velocity field to derive a global-in-time bound for the density and uniqueness of Fujita-Kato solution.
\end{Remark}


\subsection{Proof strategy} \ \
\ \   Before going into the heart of the proof of our main results,  we make a brief interpretation of the main difficulties and strategies involved in the proof.

From a mathematical point of view, the absence of the pressure term
$P$ in the system \eqref{1.12} shall cause some significant challenges in mathematical analysis. In particular,  some important methods that work effectively for  the compressible viscous system \eqref{1.11} often fail for the current  pressureless system \eqref{1.12}. For instance, for the compressible viscous system \eqref{1.11},  technique  of the so-called effective viscous flux defined as $F=\Dv u-P$ exhibiting better regularity than $\Dv u$ and $P$
considered separately (see \cite{FNP, HS2, DH1, HLX,  Lions}), and the classical perturbation theory developed by Matsumura-Nishida  \cite{MN1,MN2} and  Danchin \cite{Dan1}
are ineffective here. This may partially
explain the reason why the mathematical theory of pressureless models is poorer than the classical one.


For some existing research methods regarding the pressureless system \eqref{1.12}, we  note that the method in \cite{DMP,GTZ} heavily relied
on a smallness assumption on the initial density fluctuations, and  the framework in \cite{Xu} still requires slightly more regularity on the initial density $\rho_0$, specifically that it belongs to the multiplier space  $\mathcal{M}(\dot{B}_{p,1}^{\frac{n}{p}-1})$.
Thus, the methods  for pressureless gas model \eqref{1.12}  in \cite{DMP,GTZ,Xu} are inadequate for the current case, where the initial density does not have any small perturbation assumption or regularity conditions. Therefore, we have to  adopt a completely different approach to address the problem. The new ingredients in the present paper stems from Lemma \ref{S2lem1} and its proof, which was initially utilized by Hmidi and Keraani \cite{HK} for the two-dimensional incompressible Euler system. This concept has since been further developed for three-dimensional incompressible inhomogeneous models, as seen in \cite{ QCZ,XQF,P}. However, there are very few results in this direction for the compressible case, primarily because the divergence does not vanish. In this paper,  applying this method to the   non-standard heat equation with variable coefficient
$\rho(\partial_{t}u+u\cdot\nabla u)-\mu\Delta u-(\lambda+\mu)\nabla\Dv u=0$ with the initial data, and then  combining with  the low-high frequency decomposition, Bernstein's inequality  and the continuous argument,  we successfully obtain global \emph{a priori} estimates  of $\|u\|_{\tilde{L}_{t}^{\infty}(\dot{B}^{\frac12}_{2,1})}$ and $\|\nabla u\|_{\tilde{L}_{t}^{2}(\dot{B}^{\frac12}_{2,1})}$ when the initial density is merely bounded both above and below by positive constants, and the initial velocity is sufficiently small in critical Besov spaces (see Proposition \ref{3-1}).

On the other hand, the estimate of$\|\nabla u\|_{\tilde{L}_{t}^{2}(\dot{B}^{\frac12}_{2,1})}$ is no guarantee that the constructed velocity field
$u$ has a gradient in $L^{1}(\mathbb{R}^{+}; L^{\infty}(\mathbb{R}^{3}))$, which is crucial for solving hyperbolic-parabolic systems arising from fluid dynamics, particularly those with a convection-diffusion structure. For example, the global-in-time bound for the density and  the issue of uniqueness are closely tied to the Lipschitz control of the velocity field in most evolutionary compressible fluid mechanics models. Therefore, in our framework,  another  key challenge is to bound the quantity
 $\int_{0}^{\infty}\|\nabla u\|_{L^{\infty}}d\tau$.
To achieve it, employing the similar idea in the proof of Proposition \ref{3-1}, we further exploit some  extra time-weighted energy estimates for the velocity field (see Propositions \ref{3-2}-\ref{3-5}). Base on these time-weighted estimates, interpolation results, H\"{o}lder estimates in Lorentz spaces with respect to the time variable and the  crucial part of the argument involving the property that $\|t^{-\alpha} 1_{\mathbb{R}_{+}}\|_{L_{\frac{1}{\alpha}, \infty}}=1$ for any $\alpha>0$,  we eventually get Lipschitz control of the velocity field. 
Then we can immediately derive a global-in-time bound for the density in the  $L^{\infty}$-norm by using the energy method. Additionally,
for the uniqueness of solutions, the direct method based on stability
estimates for the system \eqref{1.12} is bound to fail owing to the hyperbolicity of the mass equation:
one loses one derivative, and one cannot afford any loss as $\rho$ is not regular
enough. Consequently, the uniqueness issue is non-trivial due to the roughness of the density and the hyperbolic nature of the continuity equation. To address this problem, we shall use the Lagrangian coordinates defined by the stream lines, which is motivated by \cite{D,D1,PZZ}. The advantage of this approach is that it allows us to convert the hyperbolic-parabolic coupled system into a parabolic system avoiding the loss of one derivative due to the hyperbolic part of system  \eqref{1.12}, which  makes the uniqueness issues more manageable compared to solving the system in Eulerian coordinates.

The rest of the paper is organized as follows. In the next section, we review some fundamental concepts related to Littlewood-Paley decomposition, Besov and Lorentz spaces, relevant product estimates, and other useful results. In Section \ref{ex}, we provide the proof of global existence as stated in Theorem \ref{1-1}. Finally, in the last section, we establish the uniqueness result presented in Theorem \ref{1-1}.

\noindent{\bf Notations.}  We assume $C$
to be a positive generic constant that may vary in different contexts throughout this paper. We denote $A\le CB$ by  $A\lesssim B$.

\section{Preliminaries}\ \
\ \  This section reviews various tools, including the Littlewood-Paley decomposition, Besov and Lorentz spaces, relevant product estimates, and useful propositions that will be referenced throughout the paper.

We start with the Littlewood-Paley decomposition. Choose a radial function
$\varphi\in\mathcal{S}(\mathbb{R}^{3})$ valued in the interval $[0,1]$ and
supported in
$\mathcal{C}=\{\xi\in\mathbb{R}^{3},\frac{3}{4}\leq|\xi|\leq\frac{8}{3}\}$ such
that
\begin{equation*}
\sum\limits_{j\in\mathbb{Z}}\varphi(2^{-j}\xi)=1\quad\mbox{for}\,\,  \xi\in
\mathbb{R}^{3}\setminus\{0\}.
\end{equation*}
For the  tempered distribution $f\in\mathcal{S}'(\mathbb{R}^{3})$,  the homogeneous frequency localization operators $\dot{\Delta}_{j}$ and $\dot{S}_{j}$ are defined by
\begin{equation*}
\dot{\Delta}_{j}f=\varphi(2^{-j}D)f\quad \mbox{for}\quad
\dot{S}_{j}f=\sum\limits_{j'\leq j-1}\dot{\Delta}_{j'}f.
\end{equation*}
With our choice of $\varphi$, one can easily verify that
\begin{equation*}
\dot{\Delta}_{q}\dot{\Delta}_{k}f=0 \quad \mbox{for}\,\,
|q-k|\geq2\quad\mbox{and}\quad
\dot{\Delta}_{q}(\dot{S}_{k-1}f\dot{\Delta}_{k}f)=0 \quad\mbox{for}\,\,
|q-k|\geq5.
\end{equation*}
We denote by $\mathcal{S}'_{h}(\mathbb{R}^{3})$ the space of tempered
distribution such that,  for any $f\in\mathcal{S}'(\mathbb{R}^{3})$,
\begin{equation*}
\lim\limits_{j\rightarrow-\infty}\dot{S}_{j}f=0
\end{equation*}
in the distributional sense. For any $f\in\mathcal{S}'(\mathbb{R}^{3})$,  the
following decomposition holds
\begin{equation*}
f=\sum\limits_{j\in\mathbb{Z}}\dot{\Delta}_{j}f,\quad f\in\mathcal{S}_{h}'(\mathbb{R}^{3}).
\end{equation*}

We first recall the definition of homogeneous
Besov spaces  (see, e.g.,  \cite{BCD, miao2}).
\begin{Definition}
Let $s\in\mathbb{R}$ and $1\leq p,r\leq +\infty$. The homogeneous Besov space
$\dot{B}_{p,r}^{s}(\mathbb{R}^{3})$ is defined by
\begin{equation*}
\dot{B}_{p,r}^{s}(\mathbb{R}^{3})=\Big{\{}f\in\mathcal{S}^{'}_{h}(\mathbb{R}^{3}):\|f\|_{\dot{B}_{p,r}^{s}(\mathbb{R}^{3})}<\infty\Big{\}},
\end{equation*}
where
\begin{equation*}
\|f\|_{\dot{B}_{p,r}^{s}(\mathbb{R}^{3})}\eqdefa\big{\|}2^{js}\|\dot{\triangle}_{j}f\|_{L^{p}}\big{\|}_{l^{r}}.
\end{equation*}
\end{Definition}

We also  introduce the following  Besov-Chemin-Lerner spaces
$\widetilde{L}^q_T(\dot{B}^{s}_{p,r})$ which is initiated in
\cite{Che-Ler}.
\begin{Definition}\label{66} Let $s\leq\frac{3}{p}$ (respectively $s\in\mathbb{R}$),
$(r,\rho,p)\in[1,\,+\infty]^3$ and $T\in(0,\,+\infty]$. We
define $\widetilde{L}^{\rho}_T(\dot B^s_{p\,r})$ as the completion of  $C([0, T]; \mathcal{S}_{h}')$  by the norm
$$
\|f\|_{\tilde{L}^\rho_{T}(\dot{B}_{p,  r}^s)}\eqdefa \Big\| 2^{js}
\|\dot{\Delta}_j f(t)\|_{L^\rho(0,T;L^{p})}\Big\|_{\ell^r}
<\infty,
$$
with the usual change if $r=\infty.$ 
\end{Definition}

Obviously, $\tilde{L}^1_{T}(\dot{B}_{p,1}^s)=L^1_{T}(\dot{B}_{p,1}^s)$, by a direct application of Minkowski's inequality, we have the following relations
between these spaces
$$L^\rho_{T}(\dot{B}_{p,r}^s)\hookrightarrow \tilde{L}^\rho_{T}(\dot{B}_{p,  r}^s),\quad  r\geq \rho,$$
$$\tilde{L}^\rho_{T}(\dot{B}_{p,  r}^s)\hookrightarrow L^\rho_{T}(\dot{B}_{p,r}^s) ,\quad  \rho\geq r.$$

Here, we also  present the definition of the Lorentz space.
\begin{Definition}\label{2ddd-1.2}
Given $f$ a measurable function on a measure space $(X, \mu)$ and $1\leq p,r\leq\infty$, we define
\begin{equation*}\label{2.1}
\tilde{\|}f\|_{L^{p,r}(\mathbb{X},\mu)}=\left\{
\begin{array}{cl}
\Big(\int_{0}^{\infty}\big(t^{\frac{1}{p}}f^{*}(t)\big)^{r}\frac{dt}{t}\Big)^{\frac{1}{r}}\quad&\text{if}\quad r<\infty,\\
\sup_{t>0}\;t^{\frac{1}{p}}f^{*}(t)\quad&\text{if}\quad r=\infty,\\
\end{array}\right.
\end{equation*}
where
$$f^{*}(t):=\inf\big\{s\geq0:\mu(\{|f|>s\})\leq t\big\}.$$
The set of all $f$ with $\tilde{\|}f\|_{L^{p,r}(\mathbb{X},\mu)}<\infty$ is called the Lorentz space with indices $p$ and $r$.
\end{Definition}
\begin{Remark}\label{2d-1.2}
Since $L^{p,p}(X,\mu)$ coincide with the Lebesgue space $L^{p}(X,\mu)$. Then, the Lorentz spaces may be endowed with the following quasi-norm,
\begin{equation*}\label{2.1}
\|f\|_{L^{p,r}(\mathbb{X},\mu)}=\left\{
\begin{array}{cl}
p^{\frac{1}{r}}\Big(\int_{0}^{\infty}\big(s\;\mu(\{|f|>s\})^{\frac{1}{p}}\big)^{r}\frac{ds}{s}\Big)^{\frac{1}{r}}\quad&\text{if}\quad r<\infty,\\
\sup_{s>0}s\;\mu(\{|f|>s\})^{\frac{1}{p}}\quad&\text{if}\quad r=\infty.
\end{array}\right.
\end{equation*}
\end{Remark}
The following classical properties of Lorentz spaces  $L^{p,q}$ may be found in \cite{Gra}.
\begin{Proposition}\label{2-1.3}
 For $1<p, p_{1}, p_{2}<\infty$ and $1 \leq r, r_{1}, r_{2} \leq \infty$, we have

(1) Interpolation. For all $1 \leq r, q \leq \infty$ and $\theta \in(0,1)$, we have
$$\big(L^{p_{1}}(\mathbb{R}_{+} ; L^{q}(\mathbb{R}^{3})); L^{p_{2}}(\mathbb{R}_{+} ; L^{q}(\mathbb{R}^{3}))\big)_{\theta, r}=L^{p, r}\big(\mathbb{R}_{+} ; L^{q}(\mathbb{R}^{3})\big),$$
where $1<p_{1}<p<p_{2}<\infty$ are such that $\frac{1}{p}=\frac{(1-\theta)}{p_{1}}+\frac{\theta}{p_{2}}$.

(2) H\"{o}lder inequality. If $\frac{1}{p}=\frac{1}{p_{1}}+\frac{1}{p_{2}}$, and $ \frac{1}{r}=\frac{1}{r_{1}}+\frac{1}{r_{2}},$ we have
$\|f g\|_{L^{p, r}} \lesssim\|f\|_{L^{p_{1}, r_{1}}}\|g\|_{L^{p_{2}, r_{2}}}.$

(3) For any $t>0$, we have $\|t^{-\alpha} 1_{\mathbb{R}_{+}}\|_{L^{\frac{1}{\alpha}, \infty}}=1$.
\end{Proposition}

The following Bernstein's inequality will be frequently used (see \cite{CH}).
\begin{Lemma}\label{Lem:Bernstein}
Let $1\le p\le q\le+\infty$. Assume that $f\in L^p(\mathbb{R}^3)$,
then for any $\gamma\in(\mathbb{N}\cup\{0\})^3$, there exist
constants $C_1$, $C_2$ independent of $f$, $j$ such that \beno
&&{\rm supp}\hat f\subseteq \{|\xi|\le A_02^{j}\}\Rightarrow
\|\partial^\gamma f\|_q\le C_12^{j{|\gamma|}+j
3(\frac{1}{p}-\frac{1}{q})}\|f\|_{p};
\\
&&{\rm supp}\hat f\subseteq \{A_12^{j}\le|\xi|\le
A_22^{j}\}\Rightarrow \|f\|_{p}\le
C_22^{-j|\gamma|}\sup_{|\beta|=|\gamma|}\|\partial^\beta f\|_p.
\eeno
\end{Lemma}

The usual product is continuous in many Besov spaces. The proof of the following lemma
 which may be found in \cite{RS} section 4.4 (
see in particular inequality (28) page 174).
\begin{Lemma}\label{p26}
For all $1\leq r,p\leq+\infty$, there exists a positive universal constant $C$  such that
$$\|fg\|_{\dot{B}^{s}_{p,r}}\leq C
\|f\|_{L^{\infty}}\|g\|_{\dot{B}^{s}_{p,r}}
+C\|g\|_{L^{\infty}}\|f\|_{\dot{B}^{s}_{p,r}}, \quad
\text{if}\quad s>0;$$
$$\|fg\|_{\dot{B}^{s_1+s_2-\frac{d}{p}}_{p,r}}\leq C
\|f\|_{\dot{B}^{s_1}_{p,r}}\|g\|_{\dot{B}^{s_2}_{p,\infty}}, \quad
\text{if}\quad s_1,s_2<\frac{3}{p},\quad \text{and}\quad
s_1+s_2>0;$$
$$\|fg\|_{\dot{B}^{s}_{p,r}}\leq C
\|f\|_{\dot{B}^{s}_{p,r}}\|g\|_{\dot{B}^{\frac{3}{p}}_{p,\infty}\cap L^{\infty}}, \quad
\text{if}\quad |s|<\frac{3}{p}.$$
\end{Lemma}

Some embedding properties  and  interpolation inequalities   about the  Besov spaces  which may be found in \cite{BCD} are in order:
\begin{Lemma}\label{equ:lemma101}
\begin{itemize}
\item\,\, For any $p\in[1,\infty]$, we have the continuous embedding
$$\dot{B}_{p,1}^{0}\hookrightarrow L^{p}\hookrightarrow \dot{B}_{p,\infty}^{0}.$$\par
\item\,\, If $s\in\mathbb{R},\quad 1\leq p_{1}\leq p_{2}\leq\infty$ and $1\leq r_{1}\leq r_{2}\leq\infty$, then $\dot{B}_{p_{1},r_{1}}^{s}\hookrightarrow
\dot{B}_{p_{2},r_{2}}^{s-3(\frac{1}{p_{1}}-\frac{1}{p_{2}})}.$\par
\item\,\, The space $\dot{B}_{p,1}^{\frac{3}{p}}$ is continuously embedded in the set of bounded continuous functions \par
(going to 0 at infinity if $p<\infty$).
\item\,\, For $1\leq p, r_{1}, r_{2}, r \leq \infty, $ $\sigma_{1}\neq \sigma_{2}$ and $\theta\in(0,1)$, then
$$\|f\|_{\dot{B}_{p,r}^{\theta\sigma_{2}+(1-\theta)\sigma_{1}}}\leq C
\|f\|^{1-\theta}_{\dot{B}_{p,r_{1}}^{\sigma_{1}}}
\|f\|^{\theta}_{\dot{B}_{p,r_{2}}^{\sigma_{2}}}.$$
\end{itemize}
\end{Lemma}

The main idea to the proof  of Theorem \ref{1-1}  is motivated by the following lemma and its proof in \cite{BCD}.
\begin{Lemma}\cite[Lemma 2.64]{BCD}\label{S2lem1}
 Let $s$ be a positive real number and $(p,r)$ be in $[1,\infty]^2.$ A constant $C_s$ exits such that
if $(u_j)_{j\in\mathbb{Z}}$ is a sequence of smooth functions where $\sum_{j\in\mathbb{Z}}u_j$ converges to $u$ in $\mathcal{S}'_{h}(\mathbb{R}^{3})$ and
$$
N_s\bigl((u_j)_{j\in\mathbb{Z}}\bigr)\eqdefa \bigl\|\bigl(\sup_{|\al|\in \{0,[s]+1\}}2^{j(s-|\al|)}\|\p^\al u_j\|_{L^p}\bigr)_{j\in\mathbb{Z}}\bigr\|_{\ell^r(\mathbb{Z})}<\infty,
$$
then $u$ is in $B^s_{p,r}$ and $\|u\|_{B^s_{p,r}}\leq C_sN_s\bigl((u_j)_{j\in\mathbb{Z}}\bigr).$
\end{Lemma}

\section{The proof of global existence in Theorem \ref{1-1}}
\label{ex}
\ \ \ \ This section is dedicated to proving the global existence of solutions to the Cauchy problem \eqref{1.12}-\eqref{1.13}. Due to the length of the proof, it has been divided into several subsections for clarity.
\subsection{Global \emph{a priori}  estimates for the velocity $u$}
\ \ \ \  Assuming that \eqref{1.4} holds, we will first establish the global \emph{a priori} estimates for the velocity $u$ by combining the energy method with the proof of Lemma \ref{S2lem1}.
\begin{Proposition}\label{3-1} Under the assumptions of \eqref{1.2} and \eqref{1.3}, let $(\rho,u)$ be a smooth enough solution of  the Cauchy problem \eqref{1.12}-\eqref{1.13}.  Then we have
\begin{equation}\label{3.1}
\|u\|_{\tilde{L}_{t}^{\infty}(\dot{B}^{\frac12}_{2,1})}+\|\nabla u\|_{\tilde{L}_{t}^{2}(\dot{B}^{\frac12}_{2,1})}\leq C\|u_{0}\|_{\dot{B}^{\frac12}_{2,1}}\quad \text{for}\quad t>0.
\end{equation}
\end{Proposition}
\noindent{\bf Proof.} To bound \eqref{3.1}, we first study the following coupled system:
\begin{equation}\label{3.2}
\left\{
\begin{aligned}
&\rho(\partial_{t}u_{j}+u\cdot\nabla u_{j})-\mu\Delta u_{j}-(\lambda+\mu)\nabla\Dv u_{j}=0,\\
&u_{j}|_{t=0}=\Delta_{j}u_0,
\end{aligned}
\right.
\end{equation}
where $\{m_{j}\}_{j\in\mathbb{Z}}$ is a sequence of smooth functions satisfying $\mathrm{supp} m_{j}\subset 2^{j}\mathcal{C}$ and $\Big\|\{2^{js}\|m_{j}\|_{L^{p}}\}_{j\in\mathbb{Z}}\Big\|_{l^{r}}<\infty$ with $0<s<\frac{d}{p}(1\leq p,r\leq\infty)$ or $s=\frac{d}{p}(r=1)$.
Then we deduce from the uniqueness of local smooth solution to the Cauchy problem \eqref{1.12}-\eqref{1.13} on $[0, T^{\ast})$ that
\begin{equation}\label{3.3}
u=\sum\limits_{j\in\mathbb{Z}}u_{j},
\end{equation}
which together the low-high frequency decomposition and Bernstein's inequality in Lemma \ref{Lem:Bernstein} immediately  yields that
\begin{equation}
\begin{split}\label{3.4}
&\|\dot{\Delta}_{j}u\|_{L_{t}^{\infty}(L^{2})}+\|\nabla\dot{\Delta}_{j}u\|_{L_{t}^{2}(L^{2})}\\
&\quad\lesssim\sum\limits_{j'>j}\Big(\|\dot{\Delta}_{j}u_{j'}\|_{L_{t}^{\infty}(L^{2})}+\|\nabla\dot{\Delta}{j}u_{j'}\|_{L_{t}^{2}(L^{2})}\Big)\\
&\qquad+2^{-j}\sum\limits_{j'\leq j}\Big(\|\nabla\dot{\Delta}_{j}u_{j'}\|_{L_{t}^{\infty}(L^{2})}+\|\nabla^{2}\dot{\Delta}_{j}u_{j'}\|_{L_{t}^{2}(L^{2})}\Big)\\
&\quad\lesssim\sum\limits_{j'>j}\Big(\|u_{j'}\|_{L_{t}^{\infty}(L^{2})}+\|\nabla u_{j'}\|_{L_{t}^{2}(L^{2})}\Big)\\
&\qquad+2^{-j}\sum\limits_{j'\leq j}\Big(\|\nabla u_{j'}\|_{L_{t}^{\infty}(L^{2})}+\|\nabla^{2}u_{j'}\|_{L_{t}^{2}(L^{2})}\Big).
\end{split}
\end{equation}
Now, we first bound the terms $\|u_{j}\|_{L_{t}^{\infty}(L^{2})}+\|\nabla u_{j}\|_{L_{t}^{2}(L^{2})}$ in the above inequality  as follows. To begin with, taking the $L^2$-scalar product of the first equation of the system \eqref{3.2} with $u_{j}$, we write
\begin{equation*}
\frac{1}{2}\frac{d}{dt}\int_{\mathbb{R}^{3}}\rho|u_{j}|^{2}dx+\mu\int_{\mathbb{R}^{3}}|\nabla u_{j}|^{2}dx+(\lambda+\mu)\int_{\mathbb{R}^{3}}|\Dv u_{j}|^{2}dx=0.
\end{equation*}
Integrating the above equation over $[0,t]$,  we have
\begin{equation*}
\frac12\|\sqrt{\rho}u_{j}\|_{L_{t}^{\infty}(L^{2})}^{2}+\mu\|\nabla u_{j}\|_{L_{t}^{2}(L^{2})}^{2}+(\lambda+\mu)\|\Dv u_{j}\|_{L_{t}^{2}(L^{2})}^{2}= \frac12\|\sqrt{\rho}\dot{\Delta}_{j}u_{0}\|_{L^{2}}^{2},
\end{equation*}
and thus there holds
\begin{equation}\label{3.5}
\|u_{j}\|_{L_{t}^{\infty}(L^{2})}+\|\nabla u_{j}\|_{L_{t}^{2}(L^{2})}+\|\Dv u_{j}\|_{L_{t}^{2}(L^{2})}\leq C\|\dot{\Delta}_{j}u_{0}\|_{L^{2}}
\leq Cd_{j}2^{-\frac{j}{2}}\|u_{0}\|_{\dot{B}^{\frac12}_{2,1}}
\end{equation}
with $\sum_{j\in\mathbb{Z}}d_{j}=1$.

For the term $\|\nabla u_{j}\|_{L_{t}^{\infty}(L^{2})}+\|\nabla^{2}u_{j}\|_{L_{t}^{2}(L^{2})}$,  taking $L^2$ inner product of the momentum equation of \eqref{3.2} with $\partial_t u_j$, and using H\"{o}lder's inequality and  $\dot{H}^{\frac{1}{2}}(\mathbb{R}^{3})\hookrightarrow L^{3}(\mathbb{R}^{3})$, we obtain
\begin{equation}
\begin{split}\label{3.6}
&\|\sqrt{\rho}\partial_{t}u_{j}\|_{L^{2}}^{2}+\frac{\mu}{2}\frac{d}{dt}\|\nabla u_{j}\|_{L^{2}}^{2}+\frac{\lambda+\mu}{2}\frac{d}{dt}\|\Dv u_{j}\|_{L^{2}}^{2}\\
&=\int_{\mathbb{R}^{3}}-\rho u\cdot\nabla u_{j}\cdot\partial_{t} u_{j}dx\\
&\leq C\|u\|_{L^{3}}\|\nabla u_{j}\|_{L^{6}}\|\sqrt{\rho}\partial_{t}u_{j}\|_{L^{2}}\\
&\leq C\|u\|_{\dot{H}^{\frac{1}{2}}}\|\nabla^{2}u_{j}\|_{L^{2}}\|\sqrt{\rho}\partial_{t}u_{j}\|_{L^{2}}.
\end{split}
\end{equation}
For second order derivatives of $u_{j}$, from \eqref{3.2}, we get
\begin{equation}\label{3.7}
-\mu\Delta u_{j}-(\lambda+\mu)\nabla\Dv u_{j}=-\rho\partial_{t}u_{j}-\rho u\cdot\nabla u_{j},
\end{equation}
which together with the  classical  regularity theory of  the Lam\'{e} system   implies that
\begin{equation*}
\begin{split}\label{3.8}
\|\nabla^{2}u_{j}\|_{L^{2}}&\leq C(\|\sqrt{\rho}\partial_{t}u_{j}\|_{L^{2}}+\|\rho u\cdot\nabla u_{j}\|_{L^{2}})\\
&\leq C(\|\sqrt{\rho}\partial_{t}u_{j}\|_{L^{2}}+\|u\|_{L^{3}}\|\nabla u_{j}\|_{L^{6}})\\
&\leq C(\|\sqrt{\rho}\partial_{t}u_{j}\|_{L^{2}}+\|u\|_{\dot{H}^{\frac{1}{2}}}\|\nabla^{2} u_{j}\|_{L^{2}}).
\end{split}
\end{equation*}
Denoting
\begin{equation}
\label{3.9}T_{1}^{\ast}\overset{\textrm{def}}{=}\sup\big\{t< T^{\ast}:\|u\|_{L_{t}^{\infty}(\dot{H}^{\frac{1}{2}})}\leq c_{1}\big\}.
\end{equation}
Then for $c_1$ in \eqref{3.9} being so small that $Cc_{1}\leq\frac{1}{2}$, thus we have for $t\leq T_{1}^{*}$ that
\begin{equation}\label{3.10}
\|\nabla^{2}u_{j}\|_{L^{2}}\leq C\|\sqrt{\rho}\partial_{t}u_{j}\|_{L^{2}}.
\end{equation}
Combining  \eqref{3.6} with \eqref{3.10} and $Cc_{1}\leq\frac{1}{2}$, we have for $t\leq T_{1}^{*}$ that
\begin{equation*}
\|\sqrt{\rho}\partial_{t}u_{j}\|_{L^{2}}^{2}+\frac{d}{dt}(\|\nabla u_{j}\|_{L^{2}}^{2}+\|\Dv u_{j}\|_{L^{2}}^{2})\leq0,
\end{equation*}
which together with  \eqref{3.10} gives
\begin{equation}
\begin{split}\label{3.11}
&\|\nabla u_{j}\|_{L_{t}^{\infty}(L^{2})}+\|\Dv u_{j}\|_{L_{t}^{\infty}(L^{2})}+\|(\nabla^{2}u_{j},\sqrt{\rho}\partial_{t}u_{j})\|_{L_{t}^{2}(L^{2})}\\
&\leq C\|\nabla u_{j}(0)\|_{L^{2}}\leq C\|\nabla\dot{\Delta}_{j}u_{0}\|_{L^{2}}\leq Cd_{j}2^{\frac{j}{2}}\|u_{0}\|_{\dot{B}^{\frac12}_{2,1}}.
\end{split}
\end{equation}
Substituting \eqref{3.5} and  \eqref{3.11} into \eqref{3.4}, we infer that
\begin{equation*}
\begin{split}\label{3.12}
&\|\dot{\Delta}_{j}u\|_{L_{t}^{\infty}(L^{2})}+\|\nabla\dot{\Delta}_{j}u\|_{L_{t}^{2}(L^{2})}\leq C d_{j}2^{-\frac{j}{2}}\|u_{0}\|_{\dot{B}^{\frac12}_{2,1}},
\end{split}
\end{equation*}
which  implies that
\begin{equation*}
\|u\|_{\tilde{L}_{t}^{\infty}(\dot{B}^{\frac12}_{2,1})}+\|\nabla u\|_{\tilde{L}_{t}^{2}(\dot{B}^{\frac12}_{2,1})}\leq C\|u_{0}\|_{\dot{B}^{\frac12}_{2,1}}\quad \text{for}\quad t\leq T_{1}^{*}.
\end{equation*}
Taking $\varepsilon_{0}$ in \eqref{1.4} so small that $C\|u_{0}\|_{\dot{B}^{\frac12}_{2,1}}\leq C\varepsilon_{0}\leq \frac{c_{1}}{2}$ for $c_{1}$ given by \eqref{3.9}. Thus, we obtain  by using a continuous argument, that $T^{\ast}_1$ determined by \eqref{3.9} equals any number smaller than  $T^{\ast}$. That is,
\begin{equation}\label{3.12-11}
\|u\|_{\tilde{L}_{t}^{\infty}(\dot{B}^{\frac12}_{2,1})}+\|\nabla u\|_{\tilde{L}_{t}^{2}(\dot{B}^{\frac12}_{2,1})}\leq C\|u_{0}\|_{\dot{B}^{\frac12}_{2,1}}\quad \text{for}\quad t\leq T^{\ast}.
\end{equation}

Finally, we will show $T^{\ast}=+\infty$.  Let $c_2$ be a small enough positive constant, which will be determined later on. Define
\begin{equation}\label{3.12-22}
T\eqdefa\max\Big\{t\in [0,T^{*}): \|u\|_{\tilde{L}_{t}^{\infty}(\dot{B}^{\frac12}_{2,1})}+\|\nabla u\|_{\tilde{L}_{t}^{2}(\dot{B}^{\frac12}_{2,1})}\leq c_2  \Big\}.
\end{equation}
We claim that  $T=+\infty$ provided that there holds \eqref{1.3}.   Indeed, taking $c_{2}=2C\varepsilon_{0}$, for $t\leq T$,  we deduce from \eqref{3.12-11} that
\begin{equation*}
\|u\|_{\tilde{L}_{t}^{\infty}(\dot{B}^{\frac12}_{2,1})}+\|\nabla u\|_{\tilde{L}_{t}^{2}(\dot{B}^{\frac12}_{2,1})}\leq C\|u_{0}\|_{\dot{B}^{\frac12}_{2,1}}\leq C\varepsilon_0= \frac{c_2}{2},
\end{equation*}
which contradicts \eqref{3.12-22}. Thus, we conclude that $T=T^{*}=\infty$, and there holds \eqref{3.1} for any $t>0$.  This completes the proof of Proposition \ref{3-1}.

\subsection{Some time weight estimates  for the velocity $u$}
\ \ \ \  Assuming that \eqref{1.4} holds, we  shall exploit  some time weight estimates  for the velocity $u$.
\begin{Proposition}\label{3-2}
Under the assumptions of Proposition \ref{3-1}, we have
\begin{equation}
\label{3.13}
\|\sqrt{t}\nabla u\|_{\tilde{L}_{t}^{\infty}(\dot{B}^{\frac12}_{2,1})}+\|\sqrt{t}\partial_{t}u\|_{\tilde{L}_{t}^{2}(\dot{B}^{\frac12}_{2,1})}
\leq C\|u_{0}\|_{\dot{B}^{\frac12}_{2,1}}\quad \text{for}\quad t>0.
\end{equation}
\end{Proposition}
\noindent{\bf Proof.} Similar to the derivation of  \eqref{3.4}, we have
\begin{equation}
\begin{split}\label{3.14}
&\|\sqrt{t}\dot{\Delta}_{j}\nabla u\|_{L_{t}^{\infty}(L^{2})}+\|\sqrt{t}\dot{\Delta}_{j}\partial_{t}u\|_{L_{t}^{2}(L^{2})}\\
&\lesssim\sum\limits_{j'>j}\Big(\|\sqrt{t}\dot{\Delta}_{j}\nabla u_{j'}\|_{L_{t}^{\infty}(L^{2})}+\|\sqrt{t}\dot{\Delta}_{j}\partial_{t} u_{j'}\|_{L_{t}^{2}(L^{2})}\Big)\\
&\quad+2^{-j}\sum\limits_{j'\leq j}\Big(\|\sqrt{t}\dot{\Delta}_{j}\nabla^{2} u_{j'}\|_{L_{t}^{\infty}(L^{2})}+\|\sqrt{t}\dot{\Delta}_{j}\nabla\partial_{t} u_{j'}\|_{L_{t}^{2}(L^{2})}\Big)\\
&\lesssim\sum\limits_{j'>j}\Big(\|\sqrt{t}\nabla u_{j'}\|_{L_{t}^{\infty}(L^{2})}+\|\sqrt{t}\partial_{t} u_{j'}\|_{L_{t}^{2}(L^{2})}\Big)\\
&\quad+2^{-j}\sum\limits_{j'\leq j}\Big(\|\sqrt{t}\nabla^{2}u_{j'}\|_{L_{t}^{\infty}(L^{2})}+\|\sqrt{t}\nabla\partial_{t}u_{j'}\|_{L_{t}^{2}(L^{2})}\Big).
\end{split}
\end{equation}
In what follows, we will bound the last two terms in the above inequality. For the terms $\|\sqrt{t}\nabla u_{j}\|_{L_{t}^{\infty}(L^{2})}$ and $\|\sqrt{t}\partial_{t}u_{j}\|_{L_{t}^{2}(L^{2})}$, it follows from  \eqref{3.6} that
\begin{equation*}
\begin{split}
&\frac{\mu}{2}\frac{d}{dt}\|\nabla u_{j}\|_{L^{2}}^{2}+\frac{\lambda+\mu}{2}\frac{d}{dt}\|\Dv u_{j}\|_{L^{2}}^{2}+\|\sqrt{\rho}\partial_{t}u_{j}\|_{L^{2}}^{2}\\
&\leq C\|u\|_{L^{\infty}}\|\nabla u_{j}\|_{L^{2}}\|\sqrt{\rho}\partial_{t}u_{j}\|_{L^{2}}\\
&\leq C\|u\|_{L^{\infty}}^{2}\|\nabla u_{j}\|^{2}_{L^{2}}+\frac{1}{2}\|\sqrt{\rho}\partial_{t}u_{j}\|^{2}_{L^{2}}.
\end{split}
\end{equation*}
Multiplying the above inequality by $t$  and using $\dot{B}^{\frac32}_{2,1}(\mathbb{R}^{3})\hookrightarrow L^{\infty}(\mathbb{R}^{3})$ in Lemma \ref{equ:lemma101} yield that
\begin{equation*}\label{3.11-11-1}
\begin{split}
&\frac{d}{dt}(\|\sqrt{t}\nabla u_{j}\|_{L^{2}}^{2}+\|\sqrt{t}\Dv u_{j}\|_{L^{2}}^{2})+\|\sqrt{t}\sqrt{\rho}\partial_{t}u_{j}\|_{L^{2}}^{2}\\
&\leq C\|\nabla u_{j}\|_{L^{2}}^{2}+C\|u\|_{\dot{B}^{\frac32}_{2,1}}^{2}\|\sqrt{t}\nabla u_{j}\|^{2}_{L^{2}}.
\end{split}
\end{equation*}
Applying Gronwall's inequality and then using \eqref{3.1} and \eqref{3.5}, we have
\begin{equation*}\label{3.11-11-2}
\begin{split}
&\|\sqrt{t}\nabla u_{j}\|_{L_{t}^{\infty}(L^{2})}^{2}+\|\sqrt{t}\Dv u_{j}\|_{L_{t}^{\infty}(L^{2})}^{2}+\|\sqrt{t}\sqrt{\rho}\partial_{t}u_{j}\|_{L_{t}^{2}(L^{2})}^{2}\\
&\leq C\|\nabla u_{j}\|_{L_{t}^{2}(L^{2})}^{2}\exp\Big(C\|u\|_{L^{2}_t(\dot{B}^{\frac32}_{2,1})}^{2}\Big)\\
&\leq Cd_{j}^{2}2^{-j}\|u_{0}\|^{2}_{\dot{B}^{\frac12}_{2,1}}.
\end{split}
\end{equation*}
This along with \eqref{3.10} implies that
\begin{equation}\label{3.15}
\|\sqrt{t}\nabla u_{j}\|_{L_{t}^{\infty}(L^{2})}+\|\sqrt{t}\Dv u_{j}\|_{L_{t}^{\infty}(L^{2})}+\|\sqrt{t}(\nabla^{2}u_{j}, \sqrt{\rho}\partial_{t}u_{j})\|_{L_{t}^{2}(L^{2})}\leq Cd_{j}2^{-\frac{j}{2}}\|u_{0}\|_{\dot{B}^{\frac12}_{2,1}}.
\end{equation}
Next, in order to bound $\|\sqrt{t}\nabla^{2}u_{j}\|_{L_{t}^{\infty}(L^{2})}$ and $\|\sqrt{t}\nabla\partial_{t}u_{j}\|_{L_{t}^{2}(L^{2})}$ in \eqref{3.14},  applying $\partial_{t}$ to  the first equation of the system \eqref{3.2} yields that
\begin{equation}\label{3.16}
\begin{split}
&\rho\partial_{t}^{2}u_{j}+\rho u\cdot\nabla\partial_{t}u_{j}-\mu\Delta\partial_{t}u_{j}-(\lambda+\mu)\nabla\Dv\partial_{t}u_{j}\\
&=-\rho_{t}D_{t}u_{j}-\rho \partial_{t}u\nabla u_{j}.
\end{split}
\end{equation}
Taking the $L^2$-scalar product of the first equation of the system \eqref{3.16} with $\partial_{t}u_{j}$, we obtain
\begin{equation}\label{3.17}
\begin{split}
&\frac{1}{2}\frac{d}{dt}\|\sqrt{\rho}\partial_{t}u_{j}\|_{L^{2}}^{2}+\mu\|\nabla\partial_{t}u_{j}\|_{L^{2}}^{2}+(\lambda+\mu)\|\Dv\partial_{t}u_{j}\|_{L^{2}}^{2}\\
&=-\int_{\mathbb{R}^{3}}\rho_{t}D_{t}u_{j}\partial_{t}u_{j}dx-\int_{\mathbb{R}^{3}}\rho \partial_{t}u\cdot\nabla u_{j}\partial_{t}u_{j}dx.
\end{split}
\end{equation}
We bound term by term above in what follows. For the term $\int_{\mathbb{R}^{3}}\rho_{t}D_{t}u_{j}\partial_{t}u_{j}dx$, we infer,  using the definition of material derivative, that
\begin{equation}\label{3.17-11-1}
\begin{split}
\int_{\mathbb{R}^{3}}\rho_{t}D_{t}u_{j}\partial_{t}u_{j}dx=\int_{\mathbb{R}^{3}}\rho_{t}|\partial_{t}u_{j}|^{2}dx+\int_{\mathbb{R}^{3}}\rho_{t} u\nabla u_{j}\cdot\partial_{t}u_{j}dx.
\end{split}
\end{equation}
For the term $\int_{\mathbb{R}^{3}}\rho_{t}|\partial_{t}u_{j}|^{2}dx$ in \eqref{3.17-11-1}, by virtue of the transport equation of \eqref{1.12}, we get, by using integration by parts, the embedding
$\dot{B}_{2,1}^{\frac{3}{2}}(\mathbb{R}^{3})\hookrightarrow L^{\infty}(\mathbb{R}^{3})$ in Lemma \ref{equ:lemma101}, H\"{o}lder's and Young's inequalities, that
\begin{equation*}\label{3.17-11-2}
\begin{split}
\int_{\mathbb{R}^{3}}\rho_{t}|\partial_{t}u_{j}|^{2}dx&\leq\int_{\mathbb{R}^{3}}\Dv(\rho u)|\partial_{t}u_{j}|^{2}dx\\
&\leq C\int_{\mathbb{R}^{3}}\rho u\cdot\nabla|\partial_{t}u_{j}|^{2}dx\\
&\leq C\|u\|_{L^{\infty}}\|\sqrt{\rho}\partial_{t}u_{j}\|_{L^{2}}\|\nabla\partial_{t}u_{j}\|_{L^{2}}\\
&\leq C\|u\|_{\dot{B}_{2,1}^{\frac{3}{2}}}^{2}\|\sqrt{\rho}\partial_{t}u_{j}\|_{L^{2}}^{2}+\frac{1}{16}\|\nabla\partial_{t}u_{j}\|^{2}_{L^{2}}
\end{split}
\end{equation*}
For the term $\int_{\mathbb{R}^{3}}\rho_{t} u\nabla u_{j}\cdot\partial_{t}u_{j}dx$ in \eqref{3.17-11-1}, by virtue of the transport equation of \eqref{1.12} and using integration by parts, we have
\begin{equation*}\label{3.17-11-3}
\begin{split}
\int_{\mathbb{R}^{3}}\rho_{t} u\nabla u_{j}&\cdot\partial_{t}u_{j}dx\leq\int_{\mathbb{R}^{3}}\Dv(\rho u)u\cdot\nabla u_{j}\partial_{t}u_{j}dx\\
&\leq\int_{\mathbb{R}^{3}}\rho u\nabla u\cdot\nabla u_{j}\partial_{t}u_{j}dx+\int_{\mathbb{R}^{3}}\rho u\cdot u\cdot\nabla^{2} u_{j}\partial_{t}u_{j}dx+\int_{\mathbb{R}^{3}}\rho u\cdot u\cdot\nabla u_{j}\cdot\nabla\partial_{t}u_{j}dx.
\end{split}
\end{equation*}
By using \eqref{3.10}, H\"{o}lder's inequality, embeddings $\dot{B}_{2,1}^{\frac{3}{2}}(\mathbb{R}^{3})\hookrightarrow L^{\infty}(\mathbb{R}^{3})$ and $\dot{B}_{2,1}^{\frac{1}{2}}(\mathbb{R}^{3})\hookrightarrow L^{3}(\mathbb{R}^{3})$ in Lemma \ref{equ:lemma101}, we obtain
\begin{equation*}
\begin{split}
\Big|\int_{\mathbb{R}^{3}}\rho u\nabla u\cdot\nabla u_{j}\partial_{t}u_{j}dx\Big|&\leq C\|u\|_{L^{\infty}}\|\nabla u\|_{L^{3}}\|\nabla u_{j}\|_{L^{6}}\|\sqrt{\rho}\partial_{t}u_{j}\|_{L^{2}}\\
&\leq C\|u\|_{\dot{B}_{2,1}^{\frac{3}{2}}}^{2}\|\nabla^{2} u_{j}\|_{L^{2}}\|\sqrt{\rho}\partial_{t}u_{j}\|_{L^{2}}\\
&\leq C\|u\|_{\dot{B}_{2,1}^{\frac{3}{2}}}^{2}\|\sqrt{\rho}\partial_{t}u_{j}\|^{2}_{L^{2}}.
\end{split}
\end{equation*}
Similarly, we can deal with the term $\int_{\mathbb{R}^{3}}\rho u\cdot u\cdot\nabla^{2} u_{j}\partial_{t}u_{j}dx$:
\begin{equation*}
\begin{split}
\Big|\int_{\mathbb{R}^{3}}\rho u\cdot u\cdot\nabla^{2} u_{j}\partial_{t}u_{j}dx\Big|&\leq C\|u\|_{L^{\infty}}\| u\|_{L^{\infty}}\|\nabla^{2} u_{j}\|_{L^{2}}\|\sqrt{\rho}\partial_{t}u_{j}\|_{L^{2}}\\
&\leq C\|u\|_{\dot{B}_{2,1}^{\frac{3}{2}}}^{2}\|\sqrt{\rho}\partial_{t}u_{j}\|^{2}_{L^{2}}.
\end{split}
\end{equation*}
Moreover, along the same line, we have
\begin{equation*}
\begin{split}
\Big|\int_{\mathbb{R}^{3}}\rho u\cdot u\cdot\nabla u_{j}\nabla\partial_{t}u_{j}dx\Big|&\leq C\|u\|_{L^{\infty}}\| u\|_{L^{3}}\|\nabla u_{j}\|_{L^{6}}\|\nabla\partial_{t}u_{j}\|_{L^{2}}\\
&\leq C\|u\|_{\dot{B}_{2,1}^{\frac{3}{2}}}\|u\|_{\dot{B}_{2,1}^{\frac{1}{2}}}\|\nabla^{2} u_{j}\|_{L^{2}}\|\nabla\partial_{t}u_{j}\|_{L^{2}}\\
&\leq C\|u\|_{\dot{B}_{2,1}^{\frac{3}{2}}}^{2}\|u\|_{\dot{B}_{2,1}^{\frac{1}{2}}}^{2}\|\sqrt{\rho}\partial_{t} u_{j}\|_{L^{2}}^{2}+\frac{1}{16}\|\nabla\partial_{t}u_{j}\|^{2}_{L^{2}}.
\end{split}
\end{equation*}
As a result, it comes out
\begin{equation}\label{3.24-A}
\begin{split}
\Big|\int_{\mathbb{R}^{3}}\rho_{t}D_{t}u_{j}\partial_{t}u_{j}dx\Big|\leq C(1+\|u\|_{\dot{B}_{2,1}^{\frac{1}{2}}}^{2})\|u\|_{\dot{B}_{2,1}^{\frac{3}{2}}}^{2}\|\sqrt{\rho}\partial_{t} u_{j}\|_{L^{2}}^{2}+\frac{1}{16}\|\nabla\partial u_{j}\|_{L^{2}}^{2}.
\end{split}
\end{equation}
For $\int_{\mathbb{R}^{3}}\rho\partial_{t}u\cdot\nabla u_{j}\partial_{t}u_{j}dx$,  using H\"{o}lder's and Young's inequalities, it follows  from \eqref{3.10} that
\begin{equation}\label{3.24-B}
\begin{split}
\Big|\displaystyle\int_{\mathbb{R}^{3}}\rho(\partial_{t}u\cdot\nabla u_{j})\partial_{t}u_{j}dx\Big|&\leq C\|\partial_{t}u\|_{L^{2}}\|\nabla u_{j}\|_{L^{3}}\|\partial_{t}u_{j}\|_{L^{6}}\\
&\leq C\|\partial_{t}u\|_{L^{2}}\|\nabla u_{j}\|_{L^{2}}^{\frac{1}{2}}\|\nabla^{2}u_{j}\|_{L^{2}}^{\frac{1}{2}}\|\partial_{t}u_{j}\|_{L^{6}}\\
&\leq C\|\partial_{t}u\|_{L^{2}}^{2}\|\nabla u_{j}\|_{L^{2}}\|\sqrt{\rho}\partial_{t}u_{j}\|_{L^{2}}+\frac{1}{16}\|\nabla\partial_{t}u_{j}\|_{L^{2}}^{2}.
\end{split}
\end{equation}
Thus, inserting \eqref{3.24-A}-\eqref{3.24-B} into \eqref{3.17} gives rise to
\begin{equation}\label{3.18}
\begin{split}
&\frac{d}{dt}\|\sqrt{\rho}\partial_{t}u_{j}\|_{L^{2}}^{2}+\|\nabla\partial_{t}u_{j}\|_{L^{2}}^{2}+\|\Dv\partial_{t}u_{j}\|_{L^{2}}^{2}\\
&\lesssim\big(1+\|u\|_{\dot{B}^{\frac12}_{2,1}}^{2}\big)\|u\|_{\dot{B}^{\frac32}_{2,1}}^{2}\|\sqrt{\rho}\partial_{t}u_{j}\|^{2}_{L^{2}}+\|\partial_{t}u\|_{L^{2}}^{2}\|\nabla u_{j}\|_{L^{2}}\|\sqrt{\rho}\partial_{t}u_{j}\|_{L^{2}}.
\end{split}
\end{equation}
Multiplying the above inequality by $t$, it follows from  \eqref{3.1} and Young's inequality that
\begin{equation*}\label{3.19}
\begin{split}
&\frac{d}{dt}\|\sqrt{t\rho}\partial_{t}u_{j}\|_{L^{2}}^{2}+\|\sqrt{t}\nabla\partial_{t}u_{j}\|_{L^{2}}^{2}+\|\sqrt{t}\Dv\partial_{t}u_{j}\|_{L^{2}}^{2}\\
&\lesssim \|\sqrt{\rho}\partial_{t}u_{j}\|_{L^{2}}^{2}+\|t^{\frac{1}{4}}\partial_{t}u\|_{L^{2}}^{2}\|\nabla u_{j}\|_{L^{2}}^{2}+(\|t^{\frac{1}{4}}\partial_{t}u\|_{L^{2}}^{2}+\|u\|_{\dot{B}^{\frac32}_{2,1}}^{2})\|\sqrt{t\rho}\partial_{t}u_{j}\|_{L^{2}}^{2}.
\end{split}
\end{equation*}
Then applying Gronwall's inequality yields
\begin{equation}\label{3.20}
\begin{split}
&\|\sqrt{t\rho}\partial_{t}u_{j}\|^{2}_{L_{t}^{\infty}(L^{2})}+\|\sqrt{t}\nabla\partial_{t}u_{j}\|^{2}_{L^{2}_{t}(L^{2})}+\|\sqrt{t}\Dv\partial_{t}u_{j}\|^{2}_{L^{2}_{t}(L^{2})}\\
&
\lesssim\big(\|\sqrt{\rho}\partial_{t}u_{j}\|_{L_{t}^{2}(L^{2})}^{2}+\|t^{\frac{1}{4}}\partial_{t}u\|_{L_{t}^{2}(L^{2})}^{2}\|\nabla u_{j}\|_{L_{t}^{\infty}(L^{2})}^{2}\big)\exp\Big(\|t^{\frac{1}{4}}\partial_{t}u\|_{L^{2}_{t}(L^{2})}^{2}+\|u\|_{L_{t}^{2}(\dot{B}^{\frac32}_{2,1})}^{2}\Big).
\end{split}
\end{equation}
In what follows, we deal with  the term $\|t^{\frac{1}{4}}\partial_{t}u_{j}\|_{L^{2}_{t}(L^{2})}$ in \eqref{3.20}.  We deduce from \eqref{3.11} and \eqref{3.15} that
\begin{equation*}\begin{split}
\|t^{\frac{1}{4}}\partial_{t}u_{j}\|_{L^{2}_{t}(L^{2})}&\leq\|\partial_{t}u_{j}\|_{L^{2}_{t}(L^{2})}^{\frac{1}{2}}\|\sqrt{t}\partial_{t}u_{j}\|_{L^{2}_{t}(L^{2})}^{\frac{1}{2}}\\
&\leq Cd_{j}\|u_{0}\|_{\dot{B}^{\frac12}_{2,1}},
\end{split}
\end{equation*}
which along with \eqref{3.3} ensures that
\begin{equation}\label{3.21}\|t^{\frac{1}{4}}\partial_{t}u\|_{L^{2}_{t}(L^{2})}\leq
\sum\limits_{j\in\mathbb{Z}}\|t^{\frac{1}{4}}\partial_{t}u_{j}\|_{L^{2}_{t}(L^{2})}\leq C\|u_{0}\|_{\dot{B}^{\frac12}_{2,1}}.
\end{equation}
By virtue of \eqref{3.1}, \eqref{3.11}  and \eqref{3.21}, we arrive at
\begin{equation}\label{3.22}\begin{split}
&\|\sqrt{t\rho}\partial_{t}u_{j}\|_{L_{t}^{\infty}(L^{2})}+\|\sqrt{t}\nabla\partial_{t}u_{j}\|_{L^{2}_{t}(L^{2})}+\|\sqrt{t}\Dv\partial_{t}u_{j}\|_{L^{2}_{t}(L^{2})}\\
&\leq Cd_{j}2^{\frac{j}{2}}\|u_{0}\|_{\dot{B}^{\frac12}_{2,1}}\exp\Big(C\|u_{0}\|_{\dot{B}^{\frac12}_{2,1}}\Big)\\
&\leq Cd_{j}2^{\frac{j}{2}}\|u_{0}\|_{\dot{B}^{\frac12}_{2,1}}.
\end{split}
\end{equation}
It then follows from  \eqref{3.10},  \eqref{3.21}  and \eqref{3.22}, that
\begin{equation}\label{3.23}
\|\sqrt{t}\nabla^{2}u_{j}\|_{L_{t}^{\infty}(L^{2})}+\|\sqrt{t}\nabla\partial_{t}u_{j}\|_{L^{2}_{t}(L^{2})}\leq Cd_{j}2^{\frac{j}{2}}\|u_{0}\|_{\dot{B}^{\frac12}_{2,1}}.
\end{equation}
Plugging \eqref{3.15} and \eqref{3.23} into \eqref{3.14},   we finally conclude that \eqref{3.13}  holds for any $t>0$.  This completes the proof of Proposition \ref{3-2}. \endproof
\begin{Proposition}\label{3-3}
Under the assumptions of Proposition \ref{3-1}, for any $t>0$, we have
\begin{equation}\label{3.23-11-1}
\|tD_{t}u_{j}\|_{L^{\infty}_{t}(L^{2})}+\|t\nabla D_{t}u_{j}\|_{L^{2}_{t}(L^{2})}\leq Cd_{j}2^{-\frac{j}{2}}\|u_{0}\|_{\dot{B}_{2,1}^{\frac{1}{2}}},
\end{equation}
and
\begin{equation}\label{3.23-11-2}
\|\sqrt{t}\nabla^{2}u\|_{L^{2}_{t}(L^{3})}\leq C\|u_{0}\|_{\dot{B}_{2,1}^{\frac{1}{2}}}.
\end{equation}
\end{Proposition}
\noindent{\bf Proof.} Multiplying \eqref{3.18} by $t^{2}$  gives rise to
\begin{equation*}
\begin{split}\label{3.26}
&\frac{d}{dt}\|t\sqrt{\rho}\partial_{t}u_{j}\|_{L^{2}}^{2}+\|t\nabla\partial_{t}u_{j}\|_{L^{2}}^{2}+\|t\Dv\partial_{t}u_{j}\|_{L^{2}}^{2}\\
&
\lesssim \|\sqrt{t\rho}\partial_{t}u_{j}\|_{L^{2}}^{2}+\|t^{\frac{1}{4}}\partial_{t}u\|_{L^{2}}^{2}\|\sqrt{t}\nabla u_{j}\|_{L^{2}}^{2}+\Big(\|t^{\frac{1}{4}}\partial_{t}u\|_{L^{2}}^{2}+\|u\|_{\dot{B}^{\frac32}_{2,1}}^{2}\Big)\|t\sqrt{\rho}\partial_{t}u_{j}\|_{L^{2}}^{2}.
\end{split}
\end{equation*}
Applying Gronwall's inequality yields that
\begin{equation*}
\begin{split}
&\|t\sqrt{\rho}\partial_{t}u_{j}\|_{L^{\infty}_{t}(L^{2})}^{2}+\|t\nabla\partial_{t}u_{j}\|_{L^{2}_{t}(L^{2})}^{2}+\|t\Dv\partial_{t}u_{j}\|_{L^{2}}^{2}\\
&\lesssim \big(\|\sqrt{t\rho}\partial_{t}u_{j}\|_{L^{2}_{t}(L^{2})}^{2}+\|t^{\frac{1}{4}}\partial_{t}u\|_{L^{2}_{t}(L^{2})}^{2}\|\sqrt{t}\nabla u_{j}\|_{L^{\infty}_{t}(L^{2})}^{2}\big)\exp\Big(\|t^{\frac{1}{4}}\partial_{t}u\|_{L^{2}_{t}(L^{2})}^{2}
+\|u\|_{L^{2}_{t}(\dot{B}^{\frac32}_{2,1})}^{2}\Big),
\end{split}
\end{equation*}
which together with \eqref{3.1}, \eqref{3.15} and \eqref{3.21} ensures that
\begin{equation}
\label{3.27}\|t\sqrt{\rho}\partial_{t}u_{j}\|_{L^{\infty}_{t}(L^{2})}+\|t\nabla\partial_{t}u_{j}\|_{L^{2}_{t}(L^{2})}+\|t\Dv\partial_{t}u_{j}\|_{L^{2}(L^{2})}\leq Cd_{j}2^{-\frac{j}{2}}\|u_{0}\|_{\dot{B}^{\frac12}_{2,1}}.\end{equation}
Whereas by applying the definition of $D_{t}$, and using H\"{o}lder's inequality, \eqref{3.1}, \eqref{3.10}, \eqref{3.13}, \eqref{3.15} and \eqref{3.27} gives rise to
\begin{equation*}
\begin{split}\label{3.28}
&\|tD_{t}u_{j}\|_{L^{\infty}_{t}(L^{2})}+\|t\nabla D_{t}u_{j}\|_{L^{2}_{t}(L^{2})}\\
&\lesssim\|t\partial_{t}u_{j}\|_{L^{\infty}_{t}(L^{2})}+\|tu\nabla u_{j}\|_{L^{\infty}_{t}(L^{2})}+\|t\nabla \partial_{t}u_{j}\|_{L^{2}_{t}(L^{2})}+\|t\nabla(u\nabla u_{j})\|_{L^{2}_{t}(L^{2})}\\
&\lesssim\|t\partial_{t}u_{j}\|_{L^{\infty}_{t}(L^{2})}+\|\sqrt{t}u\|_{L^{\infty}_{t}(L^{\infty})}\|\sqrt{t}\nabla u_{j}\|_{L^{\infty}_{t}(L^{2})}+\|t\nabla \partial_{t}u_{j}\|_{L^{2}_{t}(L^{2})}+\|\nabla u\|_{L_{t}^{2}(L^{3})}\|t\nabla u_{j}\|_{L^{\infty}_{t}(L^{6})}\\
&\quad+\|u\|_{L_{t}^{2}(L^{\infty})}\|t\nabla^2 u_{j}\|_{L^{\infty}_{t}(L^{2})}\\
&\lesssim\|t\partial_{t}u_{j}\|_{L^{\infty}_{t}(L^{2})}+\|\sqrt{t}u\|_{L^{\infty}_{t}(\dot{B}^{\frac32}_{2,1})}\|\sqrt{t}\nabla u_{j}\|_{L^{\infty}_{t}(L^{2})}+\|t\nabla \partial_{t}u_{j}\|_{L^{2}_{t}(L^{2})}+\|u\|_{L_{t}^{2}(\dot{B}^{\frac32}_{2,1})}\|t\sqrt{\rho}\partial_{t}u_{j}\|_{L^{\infty}_{t}(L^{2})}\\
&\leq Cd_{j}2^{-\frac{j}{2}}\|u_{0}\|_{\dot{B}^{\frac12}_{2,1}}.
\end{split}
\end{equation*}
On the other hand, we get, from \eqref{1.12}$_{2}$, that
\begin{equation}\label{3.23-22-1}
-\mu\Delta u-(\lambda+\mu)\nabla\Dv u=-\rho\partial_{t}u-\rho u\cdot\nabla u,
\end{equation}
which together with  the classical  regularity theory of  the Lam\'{e} system,    we infer, from  \eqref{3.1} and \eqref{3.13}, that
\begin{equation*}
\begin{split}\label{3.34}
\|\sqrt{t}\nabla^{2}u\|_{L^{2}_{t}(L^{3})}&\lesssim\|\sqrt{t}\partial_{t}u\|_{L^{2}_{t}(L^{3})}+\|u\|_{L^{2}_{t}(L^{\infty})}\|\sqrt{t}\nabla u\|_{L^{\infty}_{t}(L^{3})}\\
&\lesssim\|\sqrt{t}\partial_{t}u\|_{L^{2}_{t}(\dot{B}^{\frac12}_{2,1})} +\|u\|_{L^{2}_{t}(\dot{B}^{\frac32}_{2,1})}\|\sqrt{t}\nabla u\|_{L^{\infty}_{t}(\dot{B}^{\frac12}_{2,1})}\\
&\leq C\|u_{0}\|_{\dot{B}^{\frac12}_{2,1}}.
\end{split}
\end{equation*}
This completes the proof of Proposition \ref{3-3}.\endproof
\begin{Proposition}\label{3-4}
Under the assumptions of Proposition \ref{3-1}, for any $t>0$, we have
\begin{equation}\label{3.24}
\|tD_{t}u\|_{\tilde{L}_{t}^{\infty}(\dot{B}^{\frac12}_{2,1})}+\|t\nabla D_{t}u\|_{\tilde{L}_{t}^{2}(\dot{B}^{\frac12}_{2,1})}+\|\sqrt t\nabla u\|_{\tilde{L}_{t}^{2}(\dot{B}_{6,1}^{\frac{1}{2}})}\leq C\|u_{0}\|_{\dot{B}^{\frac12}_{2,1}}.
\end{equation}
\end{Proposition}
\noindent{\bf Proof.} Similar to the derivation of  \eqref{3.4}, we have
\begin{equation}\label{3.25}
\begin{split}
&\|t\dot{\Delta}_{j}D_{t}u\|_{L_{t}^{\infty}(L^{2})}+\|t\dot{\Delta}_{j}\nabla D_{t}u\|_{L_{t}^{2}(L^{2})}\\
&\quad\lesssim\sum\limits_{j'>j}\Big(\|t\dot{\Delta}_{j}D_{t}u_{j'}\|_{L_{t}^{\infty}(L^{2})}+\|t\dot{\Delta}_{j}\nabla D_{t}u_{j'}\|_{L_{t}^{2}(L^{2})}\Big)\\
&\qquad+2^{-j}\sum\limits_{j'\leq j}\Big(\|t\dot{\Delta}_{j}\nabla D_{t}u_{j'}\|_{L_{t}^{\infty}(L^{2})}+\|t\dot{\Delta}_{j}\nabla^{2}D_{t}u_{j'}\|_{L_{t}^{2}(L^{2})}\Big)\\
&\quad\lesssim\sum\limits_{j'>j}\Big(\|tD_{t}u_{j'}\|_{L_{t}^{\infty}(L^{2})}+\|t\nabla D_{t}u_{j'}\|_{L_{t}^{2}(L^{2})}\Big)\\
&\qquad+2^{-j}\sum\limits_{j'\leq j}\Big(\|t\nabla D_{t}u_{j'}\|_{L_{t}^{\infty}(L^{2})}+\|t\nabla^{2}D_{t}u_{j'}\|_{L_{t}^{2}(L^{2})}\Big).
\end{split}
\end{equation}
Now, we need to get the estimates for $\|t\nabla D_{t}u_{j}\|_{L_{t}^{\infty}(L^{2})}$ and $\|t\nabla^{2}D_{t}u_{j}\|_{L_{t}^{2}(L^{2})}$. Indeed, applying the operator $D_{t}$ to the first equation of the system \eqref{3.2} yields that
\begin{equation}\label{3.29}
\begin{split}
&\rho D^{2}_{t}u_{j}-\mu\Delta D_{t}u_{j}-(\lambda+\mu)\nabla\Dv D_{t}u_{j}\\
&=\rho\Dv u\cdot\frac{D}{Dt}u_{j}-\Delta u\nabla u_{j}-2\nabla u\nabla^{2}u_{j}-\nabla\Dv u\nabla u_{j}-\nabla u\nabla^{2}u_{j}-\Dv u\nabla\Dv u_{j}
\\&\eqdefa f_{j}.
\end{split}
\end{equation}
Taking the $L^2$-scalar product of  the system \eqref{3.29} with  $D^{2}_{t}u_{j}$, we obtain
\begin{equation}
\begin{split}\label{3.31}
\|\sqrt{\rho}D^{2}_{t}u_{j}\|_{L^{2}}^{2}-\int_{\mathbb{R}^{3}}\mu\Delta D_{t}u_{j}\cdot D^{2}_{t}u_{j}dx-\int_{\mathbb{R}^{3}}(\lambda+\mu)\nabla\Dv D_{t}u_{j}\cdot D^{2}_{t}u_{j}dx=\int_{\mathbb{R}^{3}}f_{j} D^{2}_{t}u_{j}dx.
\end{split}
\end{equation}
Moreover, using integration by parts and the definition of $D_{t}$ yields that
$$-\displaystyle\int_{\mathbb{R}^{3}}\Delta D_{t}u_{j}\cdot D^{2}_{t}u_{j}dx=\frac{1}{2}\frac{d}{dt}\int_{\mathbb{R}^{3}}|\nabla D_{t}u_{j}|^{2}dx+\int_{\mathbb{R}^{3}}\nabla D_{t}u_{j}\cdot\nabla(u\cdot\nabla D_{t}u_{j})dx,$$
and
$$-\displaystyle\int_{\mathbb{R}^{3}}\nabla\Dv D_{t}u_{j}\cdot D^{2}_{t}u_{j}dx=\frac{1}{2}\frac{d}{dt}\int_{\mathbb{R}^{3}}|\Dv D_{t}u_{j}|^{2}dx+\int_{\mathbb{R}^{3}}\Dv D_{t}u_{j}\cdot\Dv(u\cdot\nabla D_{t}u_{j})dx,$$
from which, multiplying \eqref{3.31} by $t^{2}$ and then integrating  over $[0,t]$, we conclude that
\begin{equation}
\begin{split}\label{3.34444}
&\|t\nabla D_{t}u_{j}\|_{L^{\infty}_{t}(L^{2})}^{2}+\|t\Dv D_{t}u_{j}\|_{L^{\infty}_{t}(L^{2})}^{2}+\|t\sqrt{\rho}D_{t}^{2}u_{j}\|_{L^{2}_{t}(L^{2})}^{2}\\
&\lesssim\|\sqrt{t}\nabla D_{t}u_{j}\|_{L^{2}_{t}(L^{2})}^{2}+\|\sqrt{t}\Dv D_{t}u_{j}\|_{L^{2}_{t}(L^{2})}^{2}+\int_{0}^{t}\tau^{2}\int_{\mathbb{R}^{3}}f_{j}D^{2}_{t}u_{j}dxd\tau\\
&\quad-\int_{0}^{t}\tau^{2}\int_{\mathbb{R}^{3}}\nabla D_{t}u_{j}\cdot\nabla(u\cdot\nabla D_{t}u_{j})dxd\tau-\int_{0}^{t}\tau^{2}\int_{\mathbb{R}^{3}}\Dv D_{t}u_{j}\cdot\Dv(u\cdot\nabla D_{t}u_{j})dxd\tau.
\end{split}
\end{equation}
In what follows, we  bound term by term  from the above inequality \eqref{3.34444}. We deduce from \eqref{3.7} that
\begin{equation*}
\begin{split}
&\|\sqrt{t}\nabla^{2}u_{j}\|_{L^{2}_{t}(L^{6})}\\
&\lesssim\|\sqrt{t}\partial_{t}u_{j}\|_{L^{2}_{t}(L^{6})}+\|\sqrt{t}u\cdot\nabla u_{j}\|_{L^{2}_{t}(L^{6})}\\
&\lesssim\|\sqrt{t}\nabla\partial_{t}u_{j}\|_{L^{2}_{t}(L^{2})}+\|u\|_{L^{2}_{t}(\dot{B}^{\frac32}_{2,1})}\|\sqrt{t}\nabla^{2}u_{j}\|_{L^{\infty}_{t}(L^{2})},
\end{split}
\end{equation*}
which together with \eqref{3.1} and \eqref{3.23} ensures that
\begin{equation}\label{3.33}
\|\sqrt{t}\nabla^{2}u_{j}\|_{L^{2}_{t}(L^{6})}\leq  Cd_{j}2^{\frac{j}{2}}\|u_{0}\|_{\dot{B}^{\frac12}_{2,1}}.
\end{equation}
By using H\"{o}lder inequality, \eqref{3.23-11-2} and \eqref{3.23}, we obtain
\begin{equation*}
\begin{split}
\|t\Delta u\cdot\nabla u_{j}\|_{L^{2}_{t}(L^{2})}&\lesssim\|\sqrt{t}\nabla^{2}u\|_{L^{2}_{t}(L^{3})}\|\sqrt{t}\nabla u_{j}\|_{L^{\infty}_{t}(L^{6})}\\
&\lesssim\|\sqrt{t}\nabla^{2}u\|_{L^{2}_{t}(L^{3})}\|\sqrt{t}\nabla^{2} u_{j}\|_{L^{\infty}_{t}(L^{2})}\\
&\leq Cd_{j}2^{\frac{j}{2}}\|u_{0}\|_{\dot{B}_{2,1}^{\frac{1}{2}}}.
\end{split}
\end{equation*}
Employing H\"{o}lder's inequality \eqref{3.13} and \eqref{3.33} yields that
\begin{equation*}
\begin{split}
\|t\nabla u\cdot\nabla^{2} u_{j}\|_{L^{2}_{t}(L^{2})}&\lesssim\|\sqrt{t}\nabla u\|_{L^{\infty}_{t}(L^{3})}\|\sqrt{t}\nabla^{2} u_{j}\|_{L^{2}_{t}(L^{6})}\\
&\lesssim\|\sqrt{t}\nabla u\|_{L^{\infty}_{t}(\dot{B}_{2,1}^{\frac{1}{2}})}\|\sqrt{t}\nabla^{2} u_{j}\|_{L^{2}_{t}(L^{6})}\\
&\leq Cd_{j}2^{\frac{j}{2}}\|u_{0}\|_{\dot{B}_{2,1}^{\frac{1}{2}}}.
\end{split}
\end{equation*}
The same estimates hold for $\|t\nabla\Dv u\cdot\nabla u_{j}\|_{L^{2}_{t}(L^{2})}$ and $\|t\Dv u\cdot\nabla\Dv u_{j}\|_{L^{2}_{t}(L^{2})}$.  For the term $\|t\rho\Dv u\cdot\frac{D}{Dt}u_{j}\|_{L^{2}_{t}(L^{2})}$,  it follows from the definition of $\frac{D}{Dt}$, H\"{o}lder's inequality, \eqref{3.1}, \eqref{3.13} and \eqref{3.23} that
\begin{equation*}
\begin{split}
\|t\rho\Dv u\cdot\frac{D}{Dt}u_{j}\|_{L^{2}_{t}(L^{2})}&\lesssim\|t\rho\Dv u\partial_{t}u_{j}\|_{L^{2}_{t}(L^{2})}+\|t\Dv u\cdot u\nabla u_{j}\|_{L^{2}_{t}(L^{2})}\\
&\lesssim\|\sqrt{t}\partial_{t}u_{j}\|_{L^{2}_{t}(L^{6})}\|\sqrt{t}\Dv u\|_{L^{\infty}_{t}(L^{3})}+\|u\|_{L^{2}_{t}(L^{\infty})}\|\sqrt{t}\nabla u\|_{L^{\infty}_{t}(L^{3})}\|\sqrt{t}\nabla u_{j}\|_{L^{\infty}_{t}(L^{6})}\\
&\lesssim\|\sqrt{t}\partial_{t}u_{j}\|_{L^{2}_{t}(L^{6})}\|\sqrt{t}\Dv u\|_{L^{\infty}_{t}(\dot{B}_{2,1}^{\frac{1}{2}})}+\|u\|_{L^{2}_{t}(\dot{B}_{2,1}^{\frac{3}{2}})}\|\sqrt{t}\nabla u\|_{L^{\infty}_{t}(\dot{B}_{2,1}^{\frac{1}{2}})}\|\sqrt{t}\nabla u_{j}\|_{L^{\infty}_{t}(L^{6})}\\
&\leq Cd_{j}2^{\frac{j}{2}}\|u_{0}\|_{\dot{B}_{2,1}^{\frac{1}{2}}}.
\end{split}
\end{equation*}
Then, we deduce that
\begin{equation}\label{3.34-11-1}
\|tf_{j}\|_{L^{2}_{t}(L^{2})}\leq Cd_{j}2^{\frac{j}{2}}\|u_{0}\|_{\dot{B}_{2,1}^{\frac{1}{2}}},
\end{equation}
and
\begin{equation*}
\begin{split}
\Big|\int_{0}^{t}\tau^{2}\int_{\mathbb{R}^{3}}f_{j}D^{2}_{t}u_{j}dxd\tau\Big|&\lesssim\|tf_{j}\|_{L^{2}_{t}(L^{2})}\|t\sqrt{\rho}D^{2}_{t}u_{j}\|_{L^{2}_{t}(L^{2})}\\
&\leq C\|tf_{j}\|_{L^{2}_{t}(L^{2})}^{2}+\frac{1}{16}\|t\sqrt{\rho}D^{2}_{t}u_{j}\|_{L^{2}_{t}(L^{2})}^{2}\\
&\leq Cd_{j}^{2}2^{j}\|u_{0}\|_{\dot{B}_{2,1}^{\frac{1}{2}}}^{2}+\frac{1}{16}\|t\sqrt{\rho}D^{2}_{t}u_{j}\|_{L^{2}_{t}(L^{2})}^{2}.
\end{split}
\end{equation*}
We finally bound $|\int_{0}^{t}\tau^{2}\int_{\mathbb{R}^{3}}\nabla D_{t}u_{j}\cdot\nabla(u\cdot\nabla D_{t}u_{j})dxd\tau|$ and $|\int_{0}^{t}\tau^{2}\int_{\mathbb{R}^{3}}\Dv D_{t}u_{j}\cdot\Dv(u\cdot\nabla D_{t}u_{j})dxd\tau|$ as follows. Applying \eqref{3.29}, we get
\begin{equation*}
-\mu\Delta D_{t}u_{j}-(\lambda+\mu)\nabla\Dv D_{t}u_{j}=-\rho D^{2}_{t}u_{j}+f_{j},
\end{equation*}
which together \eqref{3.34-11-1} implies that
\begin{equation}\label{3.34-11-2}
\begin{split}
\|t\nabla^{2}D_{t}u_{j}\|_{L_{t}^{2}(L^{2})}&\lesssim\|tf_{j}\|_{L_{t}^{2}(L^{2})}+\|t\rho D^{2}_{t}u_{j}\|_{L_{t}^{2}(L^{2})}\\
&\leq Cd_{j}2^{\frac{j}{2}}\|u_{0}\|_{\dot{B}_{2,1}^{\frac{1}{2}}}+C\|t\sqrt{\rho} D^{2}_{t}u_{j}\|_{L_{t}^{2}(L^{2})}.
\end{split}
\end{equation}
It then follows from   H\"{o}lder's and Young's inequalities, that
\begin{equation*}\label{3.43-11-1-1}
\begin{split}
&\Big|\int_{0}^{t}\tau^{2}\int_{\mathbb{R}^{3}}\nabla D_{t}u_{j}\cdot\nabla(u\cdot\nabla D_{t}u_{j})dxd\tau\Big|+\Big|\int_{0}^{t}\tau^{2}\int_{\mathbb{R}^{3}}\Dv D_{t}u_{j}\cdot\Dv(u\cdot\nabla D_{t}u_{j})dxd\tau\Big|\\
&\lesssim\int_{0}^{t}\|\nabla u\|_{L^{3}}\|\tau\nabla D_{t}u_{j}\|_{L^{3}}^{2}d\tau+\int_{0}^{t}\|u\|_{L^{\infty}}\|\tau\nabla D_{t}u_{j}\|_{L^{2}}\|\tau\nabla^{2}D_{t}u_{j}\|_{L^{2}}d\tau\\
&\lesssim\int_{0}^{t}\|\tau\nabla D_{t}u_{j}\|_{L^{2}}\|\tau\nabla^{2}D_{t}u_{j}\|_{L^{2}}(\|\nabla u\|_{L^{3}}+\|u\|_{L^{\infty}})d\tau\\
&\lesssim\|t\nabla^{2} D_{t}u_{j}\|_{L_{t}^{2}(L^{2})}\Big\|\|t\nabla D_{t}u_{j}\|_{L^{2}}\| u\|_{\dot{B}_{2,1}^{\frac{3}{2}}}\Big\|_{L^{2}_{t}}\\
&\leq\frac{1}{16}\|t\sqrt{\rho} D^{2}_{t}u_{j}\|_{L_{t}^{2}(L^{2})}^{2}+Cd_{j}^{2}2^{j}\|u_{0}\|^{2}_{\dot{B}_{2,1}^{\frac{1}{2}}}+C\int_{0}^{\tau}\|\tau\nabla D_{t}u_{j}\|_{L^{2}}^{2}\|u\|_{\dot{B}_{2,1}^{\frac{3}{2}}}^{2}d\tau.
\end{split}
\end{equation*}
As a result, it comes out
\begin{equation}\label{3.32-11-1}
\begin{split}
&\|t\nabla D_{t}u_{j}\|_{L^{\infty}_{t}(L^{2})}^{2}+\|t\Dv D_{t}u_{j}\|_{L^{\infty}_{t}(L^{2})}^{2}+\|t\sqrt{\rho}D_{t}^{2}u_{j}\|_{L^{2}_{t}(L^{2})}^{2}\\
&\lesssim d_{j}^{2}2^{j}\|u_{0}\|^{2}_{\dot{B}_{2,1}^{\frac{1}{2}}}+\|\sqrt{t}\nabla D_{t}u_{j}\|_{L^{2}_{t}(L^{2})}^{2}+\int_{0}^{\tau}\|\tau\nabla D_{t}u_{j}\|_{L^{2}}^{2}\| u\|_{\dot{B}_{2,1}^{\frac{3}{2}}}^{2}d\tau.
\end{split}
\end{equation}
Note from \eqref{3.1} and \eqref{3.23} that
\begin{equation*}
\begin{split}
\|\sqrt{t}\nabla D_{t}u_{j}\|_{L^{2}_{t}(L^{2})}&\lesssim\|\sqrt{t}\nabla \partial_{t}u_{j}\|_{L^{2}_{t}(L^{2})}+\|\sqrt{t}u\cdot\nabla u_{j}\|_{L^{2}_{t}(\dot{H}^{1})}\\
&\lesssim\|\sqrt{t}\nabla \partial_{t}u_{j}\|_{L^{2}_{t}(L^{2})}+\|u\|_{L^{2}_{t}(L^{\infty})}\|\sqrt{t}\nabla u_{j}\|_{L^{\infty}_{t}(\dot{H}^{1})}\\
&\leq Cd_{j}2^{\frac{j}{2}}\|u_{0}\|_{\dot{B}^{\frac12}_{2,1}}.
\end{split}
\end{equation*}
Then applying Gronwall's inequality to \eqref{3.34444} and using \eqref{3.1} give rise to
\begin{equation}\label{3.32-11-1}
\begin{split}
&\|t\nabla D_{t}u_{j}\|_{L^{\infty}_{t}(L^{2})}^{2}+\|t\Dv D_{t}u_{j}\|_{L^{\infty}_{t}(L^{2})}^{2}+\|t\sqrt{\rho}D_{t}^{2}u_{j}\|_{L^{2}_{t}(L^{2})}^{2}\\
&\leq Cd_{j}^{2}2^{j}\|u_{0}\|^{2}_{\dot{B}_{2,1}^{\frac{1}{2}}}exp(C\|u\|^{2}_{L^{2}_{t}(\dot{B}_{2,1}^{\frac{3}{2}})})\\
&\leq Cd_{j}^{2}2^{j}\|u_{0}\|^{2}_{\dot{B}_{2,1}^{\frac{1}{2}}},
\end{split}
\end{equation}
which together with \eqref{3.34-11-2} yields that
\begin{equation}\label{3.35-1}
\|t\nabla^{2}D_{t}u_{j}\|_{L^{2}_{t}(L^{2})}\leq Cd_{j}2^{\frac{j}{2}}\|u_{0}\|_{\dot{B}_{2,1}^{\frac{1}{2}}}.
\end{equation}
Plugging  \eqref{3.23-11-1}, \eqref{3.32-11-1} and \eqref{3.35-1} into \eqref{3.25},   we conclude that \eqref{3.24} holds for $t>0$. This completes the proof of Proposition \ref{3-4}.\endproof
\begin{Proposition}\label{3-5}
Under the assumptions of Proposition \ref{3-1}, for any $t>0$, we have
\begin{equation*}\label{3-36}
\|tu_{t}\|_{L^{\infty}_{t}(\dot{B}_{2,1}^{\frac{1}{2}})}+\|\sqrt{t}\nabla^{2} u\|_{L^{2}_{t}(L^{3})}\leq C\|u_{0}\|_{\dot{B}_{2,1}^{\frac{1}{2}}}.
\end{equation*}
\end{Proposition}
\noindent{\bf Proof.} It follows from  the law of product in Besov spaces in Lemma \ref{p26},  the definition of $D_{t}$, \eqref{3.1}, \eqref{3.13}, \eqref{3.23-22-1} and \eqref{3.24},  that
\begin{equation*}\label{3.37}
\begin{split}
&\|tu_{t}\|_{L^{\infty}_{t}(\dot{B}_{2,1}^{\frac{1}{2}})}+\|\sqrt{t}\nabla^{2} u\|_{L^{2}_{t}(L^{3})}\\
&\lesssim\|tD_{t}u\|_{L^{\infty}_{t}(\dot{B}_{2,1}^{\frac{1}{2}})}+\|tu\cdot\nabla u\|_{L^{\infty}_{t}(\dot{B}_{2,1}^{\frac{1}{2}})}+\|\sqrt{t}u_{t}\|_{L^{2}_{t}(L^{3})}+\|\sqrt{t}u\cdot\nabla u\|_{L^{2}_{t}(L^{3})}\\
&\lesssim\|tD_{t}u\|_{L^{\infty}_{t}(\dot{B}_{2,1}^{\frac{1}{2}})}+\|\sqrt{t}u\|_{L^{\infty}_{t}(\dot{B}_{2,1}^{\frac{3}{2}})}\|\sqrt{t}\nabla u\|_{L^{\infty}_{t}(\dot{B}_{2,1}^{\frac{1}{2}})}\\
&\quad+\|\sqrt{t}u_{t}\|_{L^{2}_{t}(\dot{B}_{2,1}^{\frac{1}{2}})}+\|\sqrt{t}u\|_{L^{\infty}_{t}(\dot{B}_{2,1}^{\frac{3}{2}})}\|\nabla u\|_{L^{2}_{t}(\dot{B}_{2,1}^{\frac{1}{2}})}\\
&\leq C\|u_{0}\|_{\dot{B}_{2,1}^{\frac{1}{2}}}.
\end{split}
\end{equation*}
This completes the proof of Proposition \ref{3-5}.\endproof
\subsection{The estimate  for quantity $\int_{0}^{\infty}\|\nabla u\|_{L^{\infty}}d\tau$}
\ \ \ \  Based on the all estimates mentioned above, we shall bound quantity $\int_{0}^{\infty}\|\nabla u\|_{L^{\infty}}d\tau$, which is the heart of the matter.
\begin{Proposition}\label{3-66}
Under the assumptions of Proposition \ref{3-1}, we have
\begin{equation}\label{3.52-22}\int_{0}^{\infty}\|\nabla u\|_{L^{\infty}}d\tau<\infty.\end{equation}
\end{Proposition}
\noindent{\bf Proof.} Here, we can borrow the following results from \eqref{3.15}, \eqref{3.27}, \eqref{3.10}, \eqref{3.11}, \eqref{3.23} and \eqref{3.32-11-1}
\begin{equation}\label{3.2.3}
\|\sqrt{t}\nabla^{2}u_{j}\|_{L^{2}_{t}(L^{2})}+\|t\nabla\partial_{t}u_{j}\|_{L^{2}_{t}(L^{2})}+\|t\nabla^{2}u_{j}\|_{L^{\infty}_{t}(L^{2})}\leq Cd_{j}2^{-\frac{j}{2}}\|u_{0}\|_{\dot{B}_{2,1}^{\frac{1}{2}}}
\end{equation}
and
\begin{equation}\label{3.2.4}
\|\nabla^{2}u_{j}\|_{L^{2}_{t}(L^{2})}+\|\sqrt{t}\nabla^{2}u_{j}\|_{L^{\infty}_{t}(L^{2})}+\|t\nabla\partial_{t}u_{j}\|_{L^{\infty}_{t}(L^{2})}\leq Cd_{j}2^{\frac{j}{2}}\|u_{0}\|_{\dot{B}_{2,1}^{\frac{1}{2}}}.
\end{equation}
Employing the following  interpolation of the  Lorentz space (see Proposition \ref{2-1.3}(1)),
$$
\Big(L^{2}\big(0,T;L^{2}(\mathbb{R}^{3})\big),L^{\infty}\big(0,T;L^{2}(\mathbb{R}^{3})\big)\Big)_{\frac{1}{2},1}=L^{4,1}\big(0,T;L^{2}(\mathbb{R}^{3})\big),
$$
and the estimates \eqref{3.2.3}-\eqref{3.2.4}, we infer that
\begin{equation}\label{5.2-22}
\|\sqrt{t}\nabla^{2}u_{j}\|_{L^{4,1}(0,T;L^{2}(\mathbb{R}^{3}))}\leq Cd_{j}\|(u_{0},H_{0})\|_{\dot{B}_{2,1}^{\frac{1}{2}}(\mathbb{R}^{3})}. \end{equation}
On the other hand, we also  deduce from  \eqref{3.1}, \eqref{3.13}, \eqref{3.2.3},  \eqref{3.2.4},  \eqref{3.7} with  the  classical  regularity theory of the Lam\'{e} system, H\"{o}lder's inequality and  the embeddings
$\dot{B}_{2,1}^{\frac{3}{2}}(\mathbb{R}^{3})\hookrightarrow L^{\infty}(\mathbb{R}^{3})$ in Lemma \ref{equ:lemma101} and $\dot{H}^{1}(\mathbb{R}^{3})\hookrightarrow L^{6}(\mathbb{R}^{3}),$ that
\begin{equation*}
\begin{split}
\|t\nabla^{2} u_{j}\|_{L^{2}_{t}(L^{6})}&\lesssim\|t\nabla\partial_{t}u_{j}\|_{L^{2}_{t}(L^{2})}+\|u\|_{L^{2}_{t}(L^{\infty})}\|t\nabla^{2}u_{j}\|_{L^{\infty}_{t}(L^{2})}\\
&\leq Cd_{j}2^{-\frac{j}{2}}\|u_{0}\|_{\dot{B}_{2,1}^{\frac{1}{2}}}
\end{split}
\end{equation*}
and
\begin{equation*}
\begin{split}
\|t\nabla^{2} u_{j}\|_{L^{\infty}_{t}(L^{6})}&\lesssim\|t\nabla\partial_{t}u_{j}\|_{L^{\infty}_{t}(L^{2})}+\|\sqrt{t}u\|_{L^{\infty}_{t}(L^{\infty})}\|\sqrt{t}\nabla^{2}u_{j}\|_{L^{\infty}_{t}(L^{2})}\\
&\leq Cd_{j}2^{\frac{j}{2}}\|u_{0}\|_{\dot{B}_{2,1}^{\frac{1}{2}}},
\end{split}
\end{equation*}
which together with  the interpolation of the  Lorentz space (see Proposition \ref{2-1.3}(1)): $$\big(L^{2}(0,T;L^{6}(\mathbb{R}^{3})),L^{\infty}(0,T;L^{6}(\mathbb{R}^{3}))\big)_{\frac{1}{2},1}=L^{4,1}(0,T;L^{6}(\mathbb{R}^{3})),$$
yields that
\begin{equation}\label{5.2-33}
\|t\nabla^{2}u_{j}\|_{L^{4,1}(0,T;L^{6}(\mathbb{R}^{3}))}\leq Cd_{j}\|u_{0}\|_{\dot{B}_{2,1}^{\frac{1}{2}}(\mathbb{R}^{3})}.
\end{equation}
Combining with \eqref{5.2-22} and \eqref{5.2-33} and summing for all $j\in\mathbb{Z}$, we deduce  that
\begin{equation}\label{3.2.5}
\|t\nabla^{2} u\|_{L^{4,1}_{t}(L^{6})}+\|\sqrt{t}\nabla^{2} u\|_{L^{4,1}_{t}(L^{2})}\leq C\|u_{0}\|_{\dot{B}_{2,1}^{\frac{1}{2}}}.
\end{equation}
It then follows from  Gagliardo-Nirenberg inequality, H\"{o}lder's inequality in the Lorentz spaces and \eqref{3.2.5}, that
\begin{equation*}\label{3.2.6}
\begin{split}
\int_{0}^{\infty}\|\nabla u\|_{L^{\infty}}dt&\leq C\int_{0}^{\infty}\|\nabla^{2}u\|_{L^{2}}^{\frac{1}{2}}\|\nabla^{2}u\|_{L^{6}}^{\frac{1}{2}}dt\\
&\leq C\|t^{-\frac{3}{4}}\|_{L^{\frac{4}{3},\infty}}\|t\nabla^{2} u\|_{L^{4,1}_{t}(L^{6})}^{\frac{1}{2}}\|\sqrt{t}\nabla^{2} u\|_{L^{4,1}_{t}(L^{2})}^{\frac{1}{2}}\\
&\leq C\|u_{0}\|_{\dot{B}_{2,1}^{\frac{1}{2}}}.
\end{split}
\end{equation*}
This completes the proof of Proposition \ref{3-66}.\endproof

\subsection{Bounding the density}
\ \ \ \  In this subsection, we shall deal with the density. Building on the estimate of the quantity $\int_{0}^{\infty}\|\nabla u\|_{L^{\infty}}d\tau$, we can derive a global-in-time bound for the density.
\begin{Proposition}\label{3-666}
Under the assumptions of Theorem \ref{1-1},  there holds \eqref{1.4} for $(t,x)\in \mathbb{R}^+\times\mathbb{R}^3$.
\end{Proposition}
\noindent{\bf Proof.} It follows from  \eqref{1.12}$_{1}$, for $p\geq1$,  that
\begin{equation}\label{3.56-6}
\partial_{t}(\rho^{p})+\Dv (\rho^{p} u) +(p-1)\rho \Dv u=0.
\end{equation}
Integrating \eqref{3.56-6} over $\mathbb{R}^{3}$ leads to
\begin{equation*}\label{3.1.1.1.1}
\begin{split}
\frac{d}{dt}\|\rho\|_{L^{p}}\leq\frac{p-1}{p}\|\Dv u\|_{L^{\infty}}\|\rho\|_{L^{p}}.
\end{split}
\end{equation*}
For any $t>0$, then applying Gronwall's inequality gives rise to
\begin{equation*}\label{3.2.2}
\begin{split}
\|\rho\|_{L^{p}}\lesssim \frac{p-1}{p}\|\rho_{0}\|_{L^{p}}\exp\Big(\int_{0}^{t}\|\Dv u\|_{L^{\infty}}d\tau\Big).
\end{split}
\end{equation*}
Taking $p\rightarrow\infty$, we obtain
\begin{equation*}\label{3.2.2}
\begin{split}
\|\rho\|_{L^{\infty}}\lesssim\|\rho_{0}\|_{L^{\infty}}\exp\Big(\int_{0}^{t}\|\Dv u\|_{L^{\infty}}d\tau\Big),
\end{split}
\end{equation*}
which together with  \eqref{3.52-22} implies $\|\rho\|_{L^{\infty}}\leq M<\infty$,  where $M $ only depends on $C_0$.

In what follows, we derive the lower bound of the density. Let $X$ be the flow associated to the velocity $u$, that is, the solution to
\begin{equation}\label{w2}
\frac{d}{dt}X_{u}(t,y)=u(t,X_{u}(t,y)),\quad X_{u}(0,y)=y,\quad\forall y\in\mathbb{R}^{3}.
\end{equation}
Here, \eqref{w2} describes the relation between the Eulerian coordinates $x:= X_{u}(t, y)$
and the Lagrangian coordinates $y$. Thus, passing to Lagrangian coordinates, one can define directly $X_{u}$ by
$$X_{u}(t,y)=y+\int_{0}^{t}\bar{u}(\tau,y)d\tau,$$
which is a $C^{1}$ diffeomorphism over $\mathbb{R}^{3}$.
Next, we set  $J_{u}:=det(DX_{u})$ and $A_{u}=(DX_{u})^{-1}$.
Denoting
\begin{equation}\label{w3}
\eta(t,y):=\rho(t,X_{u}(t,y)),\quad v(t,y):=u(t,X_{u}(t,y))
\end{equation}
and then using \eqref{1.12}$_{1}$, and the chain rule, we find that $\eta$ satisfies
$$(J_{v}\eta)_{t}=0.$$
 Thus, we infer that
\begin{equation}\label{w333333}
\eta=J_{v}^{-1}\rho_0.
\end{equation}
By using the fact $\partial_{t}J_{v}=\Dv uJ_{v}$,  for any $t>0$, we have
$$J_{v}=\exp\big(\int_{0}^{t}\Dv u d\tau\big),$$
which together with \eqref{3.52-22} and \eqref{w333333} implies $\rho(t,x)\geq m>0$, where $m $ only depends on $c_0$. This completes the proof of Proposition \ref{3-666}.\endproof
\section{The proof of uniqueness}
\ \ \ \ With \eqref{3.52-22} at hand, we can now proceed to prove the uniqueness of solutions. The fundamental approach is to apply a Lagrangian change of coordinates to the system \eqref{1.12}, allowing us to circumvent the hyperbolic nature of the mass equation which would otherwise result in the loss of one derivative in the stability estimates.
\subsection{Lagrangian formulation}
\ \ \ \ The goal of this subsection is to recast system \eqref{1.12} in Lagrangian variables.
Recall that in light of the estimate \eqref{3.52-22}, we know that for all $k>0$, there exists a time $T_{0}>0$ such that
\begin{equation}\label{w1}
\int_{0}^{T_{0}}\|\nabla u\|_{L^{\infty}}d\tau\leq k,
\end{equation}
where the value of $k$ will be determined in the course of the computations below.

 Based on  \eqref{w2}, we firstly list some basic properties for the Lagrangian change of variables  as follows (see \cite{D}).
\begin{Proposition}\label{w-1}
Let $u$ be a velocity field with $\nabla u\in L^{1}([0,T_{0};L^{\infty}(\mathbb{R}^{3}))$, and let $X_{u}$ be its flow, defined by \eqref{w2}. Suppose that condition \eqref{w1} is fulfilled with $c<1$. Then there exists a constant $c>0$, just depending on $c$, such that the following estimates hold true, for all times $t\in[0,T_{0}]$
$$\|\mathrm{Id}-\mathrm{adj} DX_{u}(t)\|_{L^{\infty}(\mathbb{R}^{3})}\leq C\|Du\|_{L^{1}(0,T_{0};L^{\infty}(\mathbb{R}^{3}))},$$
$$\|\mathrm{Id}-A_{u}(t)\|_{L^{\infty}(\mathbb{R}^{3})}\leq C\|Du\|_{L^{1}(0,T_{0};L^{\infty}(\mathbb{R}^{3}))},$$
$$\|J^{\pm1}_{u}(t)-1\|_{L^{\infty}(\mathbb{R}^{3})}\leq C\|Du\|_{L^{1}(0,T_{0};L^{\infty}(\mathbb{R}^{3}))}.$$
\end{Proposition}
\begin{Proposition}\label{w-2}
Let $u_{1}$ and $u_{2}$ be two vector fields satisfying \eqref{w1} and define $\delta u=u^{1}-u^{2}$.  For all times $t\in[0,T_{0}]$, we have
$$\|\mathrm{adj} DX_{1}(t)-\mathrm{adj} DX_{2}(t)\|_{L^{p}(\mathbb{R}^{3})}\leq C\int_{0}^{t}\|\nabla\delta u\|_{L^{p}(\mathbb{R}^{3})}d\tau,$$
$$\|A_{1}(t)-A_{2}(t)\|_{L^{p}(\mathbb{R}^{3})}\leq C\int_{0}^{t}\|\nabla\delta u\|_{L^{p}(\mathbb{R}^{3})}d\tau,$$
$$\|J^{\pm1}_{1}(t)-J^{\pm1}_{2}(t)\|_{L^{p}(\mathbb{R}^{3})}\leq C\int_{0}^{t}\|\nabla\delta u\|_{L^{p}(\mathbb{R}^{3})}d\tau.$$
\end{Proposition}
Let us now derive the system \eqref{1.12} in the Lagrangian coordinates. By using  \eqref{w3}
 and  the chain rule, we find that $(\eta,v)$ satisfies
\begin{align} \label{w4}
\left\{
\begin{aligned}
&(J_{v}\eta)_{t}=0,\\
&\rho_{0}\partial_{t}v-\Dv_{v}(\mu\nabla_{v}v+(\mu+\lambda)(\Dv_{v}v)\mathrm{Id})=0, \\
&v|_{t=0}=u_{0},
\end{aligned}
\right.
\end{align}
where operators $\Delta_{v}$, $\nabla_{v}$ and $\Dv_{v}$ correspond to the original operators $\Delta$, $\nabla$ and $\Dv$, respectively, after performing the change to the Lagrangian coordinates. Moreover,
\begin{equation} \label{w5}
\nabla_{v}=A_{v}^{T}\nabla_{y},\quad\Dv_{v}=A_{v}^{T}:\nabla_{y}=\Dv_{y}(A_{v}\cdot),\quad\Delta_{v}=\Dv_{y}(A_{v}A_{v}^{T}\nabla_{y}\cdot).
\end{equation}
\subsection{Estimates in Lagrangian coordinates}
\ \ \ \  We are ready to prove uniqueness part in Theorem \ref{1-1} by applying the Lagrangian approach. The key ingredient is \eqref{w1}. We can taking time $T_0$ small enough so that \eqref{w1} is satisfied with $k\ll1$. We consider two solutions $(\rho,u)$ and $(\bar{\rho},\bar{u})$ of the system \eqref{1.12}, emanating from the same initial data and denoting by $(\eta,v)$ and $(\bar{\eta},\bar{v})$ the corresponding ones in Lagrangian formulation \eqref{w4}. Denote $\delta v:=\bar{v}-v$. We see that  $\delta v$ satisfying
\begin{equation} \label{4.3}
\begin{split}
\rho_{0}\delta v_{t}-&\Dv_{v}(\mu\nabla_{v}\delta v+(\mu+\lambda)(\Dv_{v}\delta v)\mathrm{Id})\\
&=\mu(\Dv_{\bar{v}}\nabla_{\bar{v}}-\Dv_{v}\nabla_{v})\bar{v}+(\mu+\lambda)(\Dv_{\bar{v}}\Dv_{\bar{v}}\mathrm{Id}-\Dv_{v}\Dv_{v}\mathrm{Id})\bar{v}
\end{split}
\end{equation}
with
$$(\Dv_{\bar{v}}\nabla_{\bar{v}}-\Dv_{v}\nabla_{v})\bar{v}=\Dv\big((adj(DX_{\bar{v}})A_{\bar{v}}^{t}-adj(DX_{v})A_{v}^{t})\cdot\nabla\bar{v}\big)$$
and
$$(\Dv_{\bar{v}}\Dv_{\bar{v}}\mathrm{Id}-\Dv_{v}\Dv_{v}\mathrm{Id})\bar{v}=\Dv\big((adj(DX_{\bar{v}})A_{\bar{v}}^{t}-adj(DX_{v})A_{v}^{t}):\nabla\bar{v}\big).$$
In what follows, we need to explain the fact that $\delta u$ is in the energy space. Obviously,   from the original system $\eqref{1.12}_{2}$, we have
$$\rho u_{t}=-\rho u\cdot\nabla u+\mu\Delta u+(\lambda+\mu)\nabla\Dv_{u}.$$
Employing \eqref{3.1}, \eqref{3.2.5}, H\"{o}lder's inequality in the Lorentz spaces, Proposition \ref{2-1.3}(3) and the embedding
$\dot{B}_{2,1}^{\frac{1}{2}}(\mathbb{R}^{3})\hookrightarrow L^{3}(\mathbb{R}^{3}),$ we have
\begin{equation}\label{06}
\begin{split}
\|\rho u_{t}\|_{L^{\frac{4}{3},1}_{t}(L^{2})}&\lesssim\|\rho u\cdot\nabla u\|_{L^{\frac{4}{3},1}_{t}(L^{2})}+\|\nabla^{2} u\|_{L^{\frac{4}{3},1}_{t}(L^{2})}\\
&\lesssim\|u\|_{L^{\infty}_{t}(L^{3})}\|\sqrt{t}\nabla^{2}u\|_{L^{4,1}_{t}(L^{2})}\|t^{-\frac{1}{2}}\|_{L^{2,\infty}}\\
&\quad+\|\sqrt{t}\nabla^{2}u\|_{L^{4,1}(L^{2})}\|t^{-\frac{1}{2}}\|_{L^{2,\infty}}\\
&<\infty,
\end{split}
\end{equation}
which together with \eqref{1.4} implies that
\begin{equation}\label{AA}
u_{t}\in L^{\frac{4}{3},1}\big(0,T;L^{2}(\mathbb{R}^{3})\big).
\end{equation}
Thus  we may deduce from $u_{t}\in L^{\frac{4}{3},1}\big(0,T;L^{2}(\mathbb{R}^{3})\big)$ and $u\in \mathcal{C}_{b}\big(0,T;\dot{B}_{2,1}^{\frac{1}{2}}(\mathbb{R}^{3})\big)$, that
$$u(t)-u_{0}\in \mathcal{C}\big(0,T;L^{2}(\mathbb{R}^{3})\big)\cap \mathcal{C}\big(0,T;\dot{B}_{2,1}^{\frac{1}{2}}(\mathbb{R}^{3})\big),$$
which implies that  $u(t)-u_{0}\in \mathcal{C}\big(0,T; B_{2,1}^{\frac{1}{2}}(\mathbb{R}^{3})\big)$ (nonhomogeneous Besov space). Owing to the classical embedding $B_{2,1}^{\frac{1}{2}}(\mathbb{R}^{3})\hookrightarrow H^{\frac{1}{2}}(\mathbb{R}^{3})$, we get
$$u(t)-u_{0}\in \mathcal{C}\big(0,T;L^{2}(\mathbb{R}^{3})\big).$$
 Now, we claim
\begin{equation}\label{07}
\begin{split}
\nabla u \in L^{4,1}(0,T;L^{2}(\mathbb{R}^{3})).
\end{split}
\end{equation}
 Indeed, one  takes two constants $q_{0}$ and $q_{1}$ such that $1<q_{0}<\frac{4}{3}<q_{1}<\infty$ and $\frac{1}{q_{0}}+\frac{1}{q_{1}}=\frac{3}{2}$. For all $\gamma \in(0,1)$ and $i=0,1$, by using the mixed derivative theorem, we get
\begin{equation}\label{091111}\dot{W}_{2,q_{i}}^{2,1}(0,T\times\mathbb{R}^{3})\triangleq\dot{W}_{2,(q_{i},q_{i})}^{2,1}(0,T\times\mathbb{R}^{3})\hookrightarrow\dot{W}_{q_{i}}^{\gamma}\big(0,T;\dot{W}_{2}^{2-2\gamma}(\mathbb{R}^{3})\big),\end{equation} where  the space $\dot{W}_{p,(q,r)}^{2,1}(0,T\times\mathbb{R}^{3})$  is defined by
$$\dot{W}_{p,(q,r)}^{2,1}(0,T\times\mathbb{R}^{3})\triangleq\Big\{f\in\mathcal{C}\big(0,T;\dot{B}_{p,r}^{2-\frac{2}{q}}(\mathbb{R}^{3})\big);f_{t},\nabla^{2}f\in L^{q,r}\big(0,T;L^{p}(\mathbb{R}^{3})\big)\Big\}.$$
Taking $\gamma=\frac{1}{2}$, it then follows from the Sobolev embedding, that
 \begin{equation}\label{091}
\begin{split}
\dot{W}^{\frac{1}{2}}_{q_{i}}\big(0,T;\dot{W}^{1}_{2}(\mathbb{R}^{3})\big)\hookrightarrow L^{s_{i}}\big(0,T;\dot{W}^{1}_{2}(\mathbb{R}^{3})\big)\quad \text{with} \quad \frac{1}{s_{i}}=\frac{1}{q_{i}}-\frac{1}{2}.
\end{split}
\end{equation}
On the other hand, from  the proof of (\cite{DMP}, Prop.2.1), we find  that
$$\dot{W}_{2,(\frac{4}{3},1)}^{2,1}(0,T\times\mathbb{R}^{3})=\Big(\dot{W}_{2,q_{0}}^{2,1}(0,T\times\mathbb{R}^{3});\dot{W}_{2,q_{1}}^{2,1}(0,T\times\mathbb{R}^{3})\Big)_{\frac{1}{2},1}.$$
By using  \eqref{091111} and \eqref{091} with $i=0$ and $i=1$, we have
\begin{equation}\label{08}
\begin{split}
\dot{W}_{2,(\frac{4}{3},1)}^{2,1}(0,T\times\mathbb{R}^{3})\hookrightarrow \Big(L^{s_{0}}\big(0,T;\dot{W}_{2}^{1}(\mathbb{R}^{3})\big);L^{s_{1}}\big(0,T;\dot{W}_{2}^{1}(\mathbb{R}^{3})\big)\Big)_{\frac{1}{2},1}.
\end{split}
\end{equation}
We notice that, the definition of $\gamma$, $s_{i}$ and $q_{i}$ ensure that
$$\frac{1}{2}(\frac{1}{s_{0}}+\frac{1}{s_{1}})=\frac{1}{2}(\frac{1}{q_{0}}+\frac{1}{q_{1}})-\gamma=\frac{1}{2}(\frac{1}{q_{0}}+\frac{1}{q_{1}})-\frac{1}{2}=\frac{1}{4}.$$
Hence, employing Proposition \ref{2-1.3}(1), we see that the interpolation space in the right of \eqref{08} is $L^{4,1}\big(0,T;\dot{W}_{2}^{1}(\mathbb{R}^{3})\big).$
That is,
\begin{equation*}\label{05}
\begin{split}
\dot{W}_{2,(\frac{4}{3},1)}^{2,1}(0,T\times\mathbb{R}^{3})\hookrightarrow L^{4,1}\big(0,T;\dot{W}_{2}^{1}(\mathbb{R}^{3})\big),
\end{split}
\end{equation*}
which together with  \eqref{06} and \eqref{AA} implies that \eqref{07}. Then combining  with  \eqref{07} and  Proposition \ref{2-1.3}(2)-(3) yields that
$$\nabla u\in L^{2}(0,T;L^{2}(\mathbb{R}^{3})). $$
Therefore,
$$\delta u\in L^{\infty}\big(0,T;L^{2}(\mathbb{R}^{3})\big)\cap L^{2}\big(0,T;\dot{H}^{1}(\mathbb{R}^{3})\big).$$
Now, taking the $L^{2}$-scalar product of the first equation of the system \eqref{4.3} with $\delta v$ and integrating by parts delivers
\begin{equation}\label{4.8}
\begin{split}
&\frac{1}{2} \frac{d}{d t}\|\sqrt{\rho_0} \delta v\|_{L_2}^2+\mu\|\nabla_v \delta v\|_{L_2}^2+(\mu+\lambda)\|\operatorname{div}_v \delta v\|_{L_2}^2 \\
&\lesssim \|(\operatorname{adj}(D X_{\bar{v}}) A_{\bar{v}}^{\top}-\operatorname{adj}(D X_v) A_v^{\top}) \cdot \nabla \bar{v}\|_{L^2}\|\nabla \delta v\|_{L_2} .
\end{split}
\end{equation}
In what follows, we shall bound  terms  coming from the right-hand side of  \eqref{4.8}. According to \eqref{w1}, \eqref{w5},  H\"{o}lder's inequality, Propositions \ref{w-1} and \ref{w-2}, we gather that
\begin{equation*}\label{4.9}
\begin{split}
& \frac{d}{d t}\|\sqrt{\rho_0} \delta v\|_{L_2}^2+\|\nabla \delta v\|_{L_2}^2 \lesssim\|\nabla \bar{v}\|_{L_{\infty}}\|\nabla \delta v\|_{L_2}\|\int_0^t \nabla \delta v d \tau\|_{L_2}.
\end{split}
\end{equation*}
Through time integration and using H\"{o}lder's inequality, we infer that
\begin{equation*}\label{4.10}
\begin{split}
&\|\sqrt{\rho_0} \delta v\|_{L_2}^2+\int_{0}^{t}\|\nabla \delta v\|_{L_2}^2 d\tau\\
&\lesssim\int_{0}^{t}\|t^{\frac{1}{2}}\nabla\bar{v}\|_{L^{\infty}}\|\nabla\delta v\|_{L^{2}}d\tau\cdot\|\nabla\delta v\|_{L^{2}_{t}(L^{2})}\\
&\lesssim\|t^{\frac{1}{2}}\nabla\bar{v}\|_{L^{2}_{t}(L^{\infty})}\|\nabla\delta v\|_{L^{2}_{t}(L^{2})}^{2}.
\end{split}
\end{equation*}
Hence, there exists $\iota>0$, which is small enough, such that if, in addition to \eqref{w1}, we have
$$
\|t^\frac{1}{2} \nabla \bar{v}\|_{L^2_{t}(L^{\infty})}\leq \iota,
$$
which implies that $\delta v \equiv 0$ on $[0, T]$, that is, $\bar{v}=v$. Due to  $\delta \eta=(J_{\bar{v}}^{-1}-J_v^{-1}) \rho_0$, we also  get $\bar{\eta}=\eta$. In light of the above arguments, in order to get uniqueness on the whole $\mathbb{R}_{+}$, it suffices to show that our solutions satisfy not only $\nabla u \in L^1(\mathbb{R}_{+} ; L^{\infty}(\mathbb{R}^3))$, but also
$$
\int_0^{\infty} t\|\nabla u\|_{L^{\infty}}^2 d\tau<\iota.
$$
Employing Gagliardo-Nirenberg inequality, H\"{o}lder's inequality in the Lorentz spaces, Proposition \ref{2-1.3}(3) and \eqref{3.2.5}, we gather that
\begin{equation*}\label{4.11}
\begin{split}
\int_0^{\infty} t\|\nabla u\|_{L^{\infty}}^2 d\tau&\lesssim\int_0^{\infty} \|t\nabla^{2}u\|_{L^{6}}\|\sqrt{t}\nabla^{2}u\|_{L^{2}}t^{-\frac{1}{2}}d\tau\\
&\lesssim\|t\nabla^{2}u\|_{L^{4,1}_{t}(L^{6})}\|\sqrt{t}\nabla^{2}u\|_{L^{4,1}_{t}(L^{2})}\|t^{-\frac{1}{2}}\|_{L^{2,\infty}}\\
&\lesssim\|u_{0}\|_{\dot{B}_{2,1}^{\frac{1}{2}}}.
\end{split}
\end{equation*}
Taking $\varepsilon_{0}$ in \eqref{1.3} is small enough, we complete the proof of uniqueness of Theorem \ref{1-1}. \endproof

\section{Declarations}



\noindent{\bf Competing interests  }\

On behalf of all authors, the corresponding author states that there is no potential conflicts of interest with respect to the research of this article.

\noindent{\bf Authors' contributions   }\

Xiaojie Wang, Jiahong Wu and Fuyi  Xu contributed equally to this work.

\noindent{\bf Funding   }\

Wang and Xu were partially supported by the
National Natural Science Foundation of China (12326430), the  Natural Science Foundation of Shandong Province (ZR2021MA017).  Wu was partially supported by the  National Science Foundation of the United
States under DMS 2104682 and DMS 2309748.

\noindent{\bf Availability of data and materials  }\

Data and materials  sharing
not applicable to this article as no data and  materials
 were generated or analyzed during the current study.

\begin{center}

\end{center}

\end{document}